\UseRawInputEncoding
\documentclass[preprint,12pt]{elsarticle}
\usepackage{geometry}
\usepackage{hyperref}
\usepackage{amsmath,amsthm,amssymb,amsfonts}
\hypersetup{hypertex=true,
colorlinks=true,
linkcolor=blue,
anchorcolor=blue,
citecolor=blue}

\usepackage{CJK,fancybox,color,graphicx,amsmath,amsthm,amssymb,enumerate,mathrsfs}
\allowdisplaybreaks[4]

\begin{document}
\numberwithin{equation}{section}
\newcommand{\D}{\displaystyle}
\newcommand{\DF}[2]{\frac{\D#1}{\D#2}}
\newtheorem{thm}{Theorem}
\newtheorem{cor}{Corollary}
\newtheorem{lem}{Lemma}
\newtheorem{prop}{Proposition}
\newtheorem{defn}[]{Definition}
\newtheorem{rem}{Remark}
\begin{frontmatter}



\title{Mixed radial-angular integrabilities for Hausdorff type operators}


\author{Ronghui Liu$^{1,2}$}
\ead{rhliu@nwnu.edu.cn}
\author{Guanghui Lu$^{1}$}
\ead{luguanghui@nwnu.edu.cn.}

\address{$^{1}$College of Mathematics and Statistics, Northwest Normal University, Lanzhou, Gansu 730070, People's Republic
of China}
\address{$^{2}$
Key Laboratory of Computing and Stochastic Mathematics (Ministry of Education), School of Mathematics and Statistics, Hunan Normal University, Changsha, Hunan 410081, People's Republic
of China}

\cortext[cor1]{ The first author is the corresponding author.}

\begin{abstract}
 This paper is devoted to studying some mixed radial-angular integrabilities for various types of Hausdorff operators and commutators.
\end{abstract}

\begin{keyword}
Hausdorff operator,  mixed radial-angular integrability,  commutator.\\
 \emph{MSC:}~~ 42B20 \sep 42B25 \sep 42B35.
\end{keyword}
\end{frontmatter}


\section{Introduction}

The discrete Hausdorff operators play  important roles in the classical summability of number series and Fourier series \cite{HS}, which can be referred to as the Hausdorff summability method and this technique was completely developed by Hausdorff \cite{Ha} by associating it with the famous moment problem. The continuous versions of Hausdorff operators were attributed to the work of Siskakas \cite{Si} in complex analysis  and Liflyand-M\'{o}ricz \cite{LM} in the Fourier transform, more details of the background and the historical development of the Hausdorff operators, we refer the reader to look the brief overview \cite{Li} and survey article \cite{CFW}.

In the previous famous results, Hausdorff operators were considered mostly on various classical function spaces, such as Lebesgue
spaces, Morrey type spaces, Herz type spaces, etc, we mention the papers \cite{GJ,GWG,GZ,WGY} and the references therein. The main works of this paper are a variety of mixed radial-angular type spaces on which the boundedness of the
Hausdorff operators are proved.

Let $S^{n-1}$ be the unit sphere in $\mathbb R^n$, $n\geq 2$, with  Lebesgue measure $d\sigma=d\sigma(\cdot)$. For any $f\in L(\mathbb R^n)$, applying the spherical coordinate formula, we write
\begin{align*}
\|f\|_{L(\mathbb R^n)}&=\int_{0}^{\infty}\int_{S^{n-1}}|f(r\theta)|d\sigma(\theta)r^{n-1}dr\\
&=\int_{S^{n-1}}\int_{0}^{\infty}|f(r\theta)|r^{n-1}drd\sigma(\theta).
\end{align*}
Inspired by this version, from the perspective of radial and angular integrabilities, we naturally consider the case of Lebesgue norms with different integrabilities in the radial and angular directions, namely, for any $1<p,\tilde{p}<\infty$,
\begin{align*}
\|f\|_{L^{p}_{\rm rad}L^{\tilde{p}}_{\rm ang}(\mathbb R^n)}
:=&\bigg(\int_{0}^{\infty}\Big(\int_{S^{n-1}}|f(r\theta)|^{\tilde{p}}d\sigma(\theta)\Big)^{p/\tilde{p}}r^{n-1}dr\bigg)^{1/p}.\\
&\neq\bigg(\int_{S^{n-1}}\Big(\int_{0}^{\infty}|f(r\theta)|^{p}r^{n-1}dr\Big)^{\tilde{p}/p}d\sigma(\theta)\bigg)^{1/\tilde{p}}=:\|f\|_{L^{\tilde{p}}_{\rm ang}L^{p}_{\rm rad}(\mathbb R^n)}.
\end{align*}
and when $p=\infty$ or $\tilde{p}=\infty$, we just need to make the usual modifications in the
above definition.

 Similar to the ideas above, the first author et al. introduced  various mixed radial-angular type  function spaces as follows:

\noindent\textbf{Definition 1.1 }(\cite{LT}) For any $1\leq p,\tilde{p}<\infty$, the weighted mixed radial-angular spaces  ${L^{p}_{rad}L^{\tilde{p}}_{ang}(\mathbb R^n, \omega)}$ are defined by
$${{L^{p}_{rad}L^{\tilde{p}}_{ang}(\mathbb R^n, \omega)}}=\big\{f:\,\|f\|_{L^{p}_{rad}L^{\tilde{p}}_{ang}(\mathbb R^n, \omega)} <\infty\big\},$$
 where
 $$\|f\|_{L^{p}_{rad}L^{\tilde{p}}_{ang}(\mathbb R^n, \omega)}=\bigg(\int_{0}^{\infty}\|f(r\cdot)\|^p_{L^{\tilde{p}}(S^{n-1})}r^{n-1}\omega(r)dr\bigg)^{{1}/{p}}.$$
For the weighted mixed radial-angular spaces $\|f\|_{L^{\tilde{p}}_{\rm ang}L^{p}_{\rm rad}(\mathbb R^n, \omega)}$, we only need to exchange the integration order in
radial and angular.

\noindent\textbf{Definition 1.2 }(\cite{LTW1})
Let $1<p_1, p_2<\infty$. A function $f\in L^{p_2}_{\rm rad,{loc}}L^{p_1}_{\rm ang}((0,\infty)$ $\times{S}^{n-1})$ is said to belong to the  mixed radial-angular  central Morrey spaces ${\mathcal{B}L^{p_2}_{\rm rad}L^{p_1}_{\rm ang}(\mathbb R^n)}$, if
\begin{small}
$$\|f\|_{{\mathcal{B}L^{p_2}_{\rm rad}L^{p_1}_{\rm ang}(\mathbb R^n)}}=\sup\limits_{r>0}\Big(\frac{1}{\nu_nr^{n}}\int_{0}^{r}\Big(\int_{{S}^{n-1}}|f(\rho\theta)|^{p_1}d\sigma(\theta)\Big)^{p_2/p_1}\rho^{n-1}d\rho\Big)^{1/p_2}<\infty,$$
\end{small}
where $f\in L^{p_2}_{\rm rad,{loc}}L^{p_1}_{\rm ang}((0,\infty)\times{S}^{n-1})$ means that the function $f$ defined on $(0,r_0)\times{S}^{n-1}\subset(0,\infty)$$\times{S}^{n-1}$ satisfies
$\|f\|_{L^{p_2}_{\rm rad}L^{p_1}_{\rm ang}((0,r_0)\times{S}^{n-1})}<\infty$ for any $r_0\in(0,\infty)$.

\noindent\textbf{Definition 1.3 }(\cite{LTW2})
Let $1<p_1, p_2<\infty$. A function $f\in L^{p_2}_{\rm  rad,{loc}}L^{p_1}_{\rm ang}((0,\infty)\times{S}^{n-1})$ is said to belong to the mixed radial-angular  CMO  spaces ${CMOL^{p_2}_{\rm rad}L^{p_1}_{\rm ang}(\mathbb R^n)}$, if
\begin{small}
\begin{equation*}
\|f\|_{{CMOL^{p_2}_{\rm rad}L^{p_1}_{\rm ang}(\mathbb R^n)}}=\sup\limits_{r>0}\Big(\frac{1}{\nu_nr^n}\int_{0}^{r}\Big(\int_{{S}^{n-1}}|f(\rho\theta)-f_{B}|^{p_1}d\sigma(\theta)\Big)^{p_2/p_1}\rho^{n-1}d\rho\Big)^{1/p_2}<\infty,
\end{equation*}
\end{small}

\noindent\textbf{Definition 1.4 }(\cite{LTW1}) Let $\alpha\in\mathbb R$, $0<p_1,p_2, q<\infty$. The mixed radial-angular homogeneous Herz spaces $\dot{K}^{\alpha,q}_{L^{p_2}_{\rm rad}L^{p_1}_{\rm ang}}(\mathbb R^{n})$
are defined by
$$\dot{K}^{\alpha,q}_{L^{p_2}_{\rm rad}L^{p_1}_{\rm ang}}(\mathbb R^{n})=\Big\{f\in L^{p_2}_{\rm rad,{loc}}L^{p_1}_{\rm ang}((0,\infty)\times{S}^{n-1}) \setminus{\{0\}}): \|f\|_{\dot{K}^{\alpha,q}_{L^{p_2}_{\rm rad}L^{p_1}_{\rm ang}}(\mathbb R^{n})}<\infty\Big\},$$
where
$$\|f\|_{{\dot{K}}^{\alpha,q}_{L^{p_2}_{\rm rad}L^{p_1}_{\rm ang}}(\mathbb R^{n})}=\Big\{\sum\limits_{k=-\infty}^{\infty}2^{k\alpha q} \|f\chi_k\|^q_{L^{p_2}_{\rm rad}L^{p_1}_{\rm ang}(\mathbb R^{n})}\Big\}^{1/q}.$$

In reality, the idea of distinguishing radial and angular directions is not new. In recent
years, the mixed radial-angular integrabilites of some classical operators has attracted the attention
of several authors \cite{LF,LLW,LLW2,LW,LRW,LTW}. The main purpose in this paper is to consider the  behaviour of  Hausdorff operators and commutators on these mixed radial-angular type function spaces, respectively.

In what follows, we use the symbol $A\lesssim B$ to denote that there exists a positive constant $C$ such that $A\leq CB$. The notation  $A\simeq B$ denotes $A$ is equivalent to $B$, that means there exist two positive constants $C_1$ and $C_2$ such that $C_1A\leq B\leq C_2A$.

The rest of this paper is organized as follows. In Section 2, we will establish various mixed radial-angular integrabilities for  Hausdorff type operators.
The results of fractional Hausdorff operators will be obtained in Section 3. The cases of  multilinear  Hausdorff operators will be established in Section 4. Finally, we will give some mixed radial-angular integrabilities for commutators of  Hausdorff operators on generalized central Morrey spaces. We would like to remark that the main ideas of our proofs
are taken from \cite{CFL,CFZ,GJ,GWG,WGY}.

\section{Mixed radial-angular integrabilities for  Hausdorff operators }

The Hausdorff operator on the Euclidean $\mathbb R^{n}$ is defined by

$$H_{\Phi}f(x)=\int_{\mathbb{R}^n}\frac{\Phi(x|y|^{-1})}{|y|^{n}}f(y)dy.$$

\noindent\textbf{Theorem 1.1} Suppose that  $\Phi$ is a radial function. For $1\leq p,\tilde{p_1},\tilde{p_2}<\infty$, we have

$$
\left\|H_{\Phi}f\right\|_{{L^{p}_{rad}L^{\tilde{p}_1}_{ang}}\left(\mathbb R^n\right)}\leq{\omega_{n-1}}^{\frac{1}{\tilde{p_1}}+\frac{1}{\tilde{p_2}'}}\int_{0}^{\infty}\frac{|\Phi(t)|}{t}t^{\frac{n}{p}}dt\left\|f\right\|_{{L^{p}_{rad}L^{\tilde{p}_2}_{ang}}\left((\mathbb R^n\right)}.
$$
\begin{proof} By polar coordinate and change of variables, we have
$$H_{\Phi}f(x)=\int_{\mathbb{R}^n}\frac{\Phi(x|y|^{-1})}{|y|^{n}}f(y)dy=\int_{0}^{\infty}\frac{\Phi(t)}{t}\int_{S^{n-1}}f(|x|t^{-1}y')d\sigma(y')dt.$$
Thus, the Minkowski inequality and  H\"{o}lder's inequality deduce that
\begin{align*}
\left\|H_{\Phi}f\right\|_{L^{p}_{rad}L^{\tilde{p}_1}_{ang}\left(\mathbb R^n \right)}\leq &\int_{0}^{\infty}\frac{|\Phi(t)|}{t}\left\|\int_{S^{n-1}}f(|\cdot|^{-1}ty')d\sigma(y')\right\|_{L^{p}_{rad}L^{\tilde{p}_1}_{ang}\left(\mathbb R^n\right)}dt \\
=&{\omega_{n-1}}^{\frac{1}{\tilde{p_1}}}\int_{0}^{\infty}\frac{|\Phi(t)|}{t}t^{\frac{n}{p}}dt\left\|f\right\|_{L^{p}_{rad}L^{1}_{ang}\left(\mathbb R^n\right)} \\
=&{\omega_{n-1}}^{\frac{1}{\tilde{p_1}}+\frac{1}{\tilde{p_2}'}}\int_{0}^{\infty}\frac{|\Phi(t)|}{t}t^{\frac{n}{p}}dt\left\|f\right\|_{L^{p}_{rad}L^{\tilde{p}_2}_{ang}\left(\mathbb R^n\right)}.
\end{align*}
Therefore, Theorem 1.1 is finished.
\end{proof}
The Hausdorff operator with rough kernel is defined by

$$H_{\Phi,\Omega}f(x)=\int_{\mathbb{R}^n}\frac{\Phi(x|y|^{-1})}{|y|^{n}}f(y)\Omega(y')dy,$$
where $\Omega(y')$ is an integrable function defined on the unit sphere $S^{n-1}$.

\noindent\textbf{Theorem 1.2} Suppose that  $\Phi$ is a radial function. For $1\leq p,\tilde{p_1},\tilde{p_2}<\infty$, $\Omega\in L^{\tilde{p}'_2}(S^{n-1})$,  we have

$$
\left\|H_{\Phi}f\right\|_{{L^{p}_{rad}L^{\tilde{p}_1}_{ang}}\left((\mathbb R^n\right)}\leq{\omega_{n-1}}^{\frac{1}{\tilde{p_1}}}\|\Omega\|_{L^{\tilde{p}'_2}(S^{n-1})}\int_{0}^{\infty}\frac{|\Phi(t)|}{t}t^{\frac{n}{p}}dt\left\|f\right\|_{{L^{p}_{rad}L^{\tilde{p}_2}_{ang}}\left(\mathbb R^n\right)}.
$$
\begin{proof} Note that
$$H_{\Phi,\Omega}f(x)=\int_{\mathbb{R}^n}\frac{\Phi(x|y|^{-1})}{|y|^{n}}f(y)\Omega(y')dy=\int_{0}^{\infty}\frac{\Phi(t)}{t}\int_{S^{n-1}}f(|x|t^{-1}y')\Omega(y')d\sigma(y')dt.$$
Thus, by the Minkowski inequality and  H\"{o}lder's inequality, we have
\begin{align*}
\left\|H_{\Phi}f\right\|_{L^{p}_{rad}L^{\tilde{p}_1}_{ang}\left(\mathbb R^n \right)}\leq &\int_{0}^{\infty}\frac{|\Phi(t)|}{t}\left\|\int_{S^{n-1}}f(|\cdot|t^{-1}y')\Omega(y')d\sigma(y')\right\|_{L^{p}_{rad}L^{\tilde{p}_1}_{ang}\left(\mathbb R^n\right)}dt \\
=&{\omega_{n-1}}^{\frac{1}{\tilde{p_1}}}\|\Omega\|_{L^{\tilde{p}'_2}(S^{n-1})}\int_{0}^{\infty}\frac{|\Phi(t)|}{t}t^{\frac{n}{p}}dt\left\|f\right\|_{L^{p}_{rad}L^{\tilde{p}_2}_{ang}\left(\mathbb R^n\right)} \\
=&{\omega_{n-1}}^{\frac{1}{\tilde{p_1}}}\|\Omega\|_{L^{\tilde{p}'_2}(S^{n-1})}\int_{0}^{\infty}\frac{|\Phi(t)|}{t}t^{\frac{n}{p}}dt\left\|f\right\|_{L^{p}_{rad}L^{\tilde{p}_2}_{ang}\left(\mathbb R^n\right)}.
\end{align*}
Therefore, Theorem 1.2 is finished.
\end{proof}

For a locally integrable function $F$ on $\mathbb{R}^{nm}$, Chen et al. \cite{CFZ} defined the following Hausdorff operator,

$$H_{\Phi}F(x)=\int_{\mathbb{R}^{nm}}\frac{\Phi(x|y|^{-1})}{|y|^{nm}}F(y)dy.$$
\noindent\textbf{Theorem 1.3} Suppose that  $\Phi$ is a radial function. Let $1\leq p,\tilde{p_1},\tilde{p_2}<\infty$ and $\beta=n(m-1)$, then we have

$$
\left\|T_{\Phi}F\right\|_{{L^{p}_{rad}L^{\tilde{p}_1}_{ang}}\left(\mathbb R^{n},  r^{\beta}\right)}\leq{\omega_{n-1}}^{\frac{1}{\tilde{p_1}}}{\omega_{nm-1}}^{\frac{1}{\tilde{p_2}'}}\int_{0}^{\infty}\frac{|\Phi(t)|}{t}t^{\frac{nm}{p}}dt\left\|F\right\|_{{L^{p}_{rad}L^{\tilde{p}_2}_{ang}}\left(\mathbb R^{nm}\right)}.
$$
\begin{proof} Similarly, we have
$$H_{\Phi}F(x)=\int_{\mathbb{R}^{nm}}\frac{\Phi(x|y|^{-1})}{|y|^{nm}}F(y)dy=\int_{0}^{\infty}\frac{\Phi(t)}{t}\int_{S^{nm-1}}F(|x|t^{-1}y')d\sigma(y')dt,$$
where $S^{nm-1}$ is the unit sphere in $\mathbb{R}^{n}$.

By the Minkowski inequality and  H\"{o}lder's inequality, we derive that
\begin{align*}
\left\|H_{\Phi}F\right\|_{L^{p}_{rad}L^{\tilde{p}_1}_{ang}\left(\mathbb R^n , r^{\beta}\right)}\leq &\int_{0}^{\infty}\frac{|\Phi(t)|}{t}\left\|\int_{S^{nm-1}}F(|\cdot|t^{-1}y')d\sigma(y')\right\|_{L^{p}_{rad}L^{\tilde{p}_1}_{ang}\left(\mathbb R^n,r^{\beta}\right)}dt \\
=&{\omega_{n-1}}^{\frac{1}{\tilde{p_1}}}\int_{0}^{\infty}\frac{|\Phi(t)|}{t}dt\Bigg(\int_{0}^{\infty}\Big(\int_{S^{nm-1}}|F(\frac{r}{t}y')|d\sigma(y')\Big)^{\frac{p}{\tilde{p_1}}}r^{nm-1}dr\Bigg)^{\frac{1}{p}}\\
=&{\omega_{n-1}}^{\frac{1}{\tilde{p_1}}}{\omega_{nm-1}}^{\frac{1}{\tilde{p_2}'}}\int_{0}^{\infty}\frac{|\Phi(t)|}{t}t^{\frac{nm}{p}}dt\left\|F\right\|_{L^{p}_{rad}L^{\tilde{p}_2}_{ang}\left(\mathbb R^{nm}\right)}.
\end{align*}
Therefore, Theorem 1.3 is finished.
\end{proof}

In Theorem 1.3, if we take
$$F(x_1,x_2,\ldots,x_m)=f_1(x_1)f_2(x_2)\ldots f_m(x_m),$$
then $H_{\Phi}F$ deduces to an $m$-linear operator
$$H_{\Phi}(f_1,f_2,\ldots,f_m)(x)=\int_{\mathbb{R}^{nm}}\frac{\Phi(x|y|^{-1})}{|y|^{nm}}\prod_{i=1}^{m}f_i(y_i)dy.$$

\noindent\textbf{Theorem 1.4} Suppose that  $\Phi$ is a radial function. Let $1\leq p,\tilde{p}<\infty$, $\frac{1}{p}=\sum_{i=1}^{m}\frac{1}{p_i}$ and $\frac{1}{\tilde{p}}=\sum_{i=1}^{m}\frac{1}{\tilde{p}_i}$, then we have

$$
\left\|H_{\Phi}(f_1,f_2,\ldots,f_m)\right\|_{{L^{p}_{rad}L^{\tilde{p}}_{ang}}\left(\mathbb R^{n}\right)}\leq\int_{\mathbb{R}^{nm}}\frac{\Phi(|y|^{-1})}{|y|^{nm}\prod_{i=1}^{m}|y_i|^{n/p_i}}\prod_{i=1}^{m}\left\|f_i\right\|_{{L^{p_i}_{rad}L^{\tilde{p}_i}_{ang}}\left(\mathbb R^{n}\right)}.
$$
\begin{proof} Lemma 2.1 in \cite{CFZ} implies that the norm of the operator of $H_{\Phi}$ is equal to the norm of the operator of $H_{\Phi}$ restricts to the set of nonnegative radial functions. Therefore, it suffices to show that the case of the function $f_i$ are nonnegative radial functions.

 By the definition and the change of variables, we  get
$$H_{\Phi}(f_1,f_2,\ldots,f_m)(x)=\int_{\mathbb{R}^{nm}}\frac{\Phi(x|y|^{-1})}{|y|^{nm}}\prod_{i=1}^{m}f_i(y_i)dy=\int_{\mathbb{R}^{nm}}\frac{\Phi(|y|^{-1})}{|y|^{nm}}\prod_{i=1}^{m}f_i(|x|y_i)dy.$$
Thus, by the Minkowski inequality and  H\"{o}lder's inequality, we have
\begin{align*}
\left\|H_{\Phi}(f_1,f_2,\ldots,f_m)\right\|_{L^{p}_{rad}L^{\tilde{p}}_{ang}\left(\mathbb R^n \right)}\leq &\int_{\mathbb{R}^{nm}}\frac{\Phi(|y|^{-1})}{|y|^{nm}}\left\|\prod_{i=1}^{m}f_i(y_i|\cdot|)\right\|_{L^{p}_{rad}L^{\tilde{p}}_{ang}\left(\mathbb R^n\right)}dy \\
=&\int_{\mathbb{R}^{nm}}\frac{\Phi(|y|^{-1})}{|y|^{nm}}\prod_{i=1}^{m}\left\|f_i(y_i|\cdot|)\right\|_{L^{p_i}_{rad}L^{\tilde{p}_i}_{ang}\left(\mathbb R^n\right)}dy\\
=&\int_{\mathbb{R}^{nm}}\frac{\Phi(|y|^{-1})}{|y|^{nm}\prod_{i=1}^{m}|y_i|^{n/p_i}}\prod_{i=1}^{m}\left\|f_i\right\|_{L^{p_i}_{rad}L^{\tilde{p}_i}_{ang}\left(\mathbb R^n\right)}.
\end{align*}
Therefore, Theorem 1.4 is finished.
\end{proof}

Chen et al. \cite{CFZ} defined another form of  Hausdorff operator in high dimension as follows:

$$\mathcal{H}_{\Phi}f(x)=\int_{\mathbb{R}^n}\frac{\Phi(y)}{|y|^{n}}f(\frac{x}{|y|})dy.$$

\noindent\textbf{Theorem 1.5} Suppose that  $f$ is a radial function. For $\alpha\in\mathbb R, 1\leq p_1, p_2, \tilde{p_1}, q<\infty$,  we have

$$
\left\|\mathcal{H}_{\Phi}f\right\|_{{\dot{K}}^{\alpha,q}_{L^{p_2}_{\rm rad}L^{p_1}_{\rm ang}}\left(\mathbb R^n\right)}\lesssim\int_{\mathbb R^{n}}\frac{|\Phi(y)|}{|y|^{n}}|y|^{\frac{n}{p_2}+\alpha}dy\left\|f\right\|_{{\dot{K}}^{\alpha,q}_{L^{p_2}_{\rm rad}L^{\tilde{p_1}}_{\rm ang}}\left(\mathbb R^n\right)}.
$$

\begin{proof} By the definition and  Minkowski inequality, we have
\begin{align*}
\left\|\mathcal{H}_{\Phi}f\right\|^q_{{\dot{K}}^{\alpha,q}_{L^{p_2}_{\rm rad}L^{p_1}_{\rm ang}}\left(\mathbb R^n\right)}&=\sum\limits_{k=-\infty}^{\infty}2^{k\alpha q} \|\mathcal{H}_{\Phi}f\chi_k\|^q_{L^{p_2}_{\rm rad}L^{p_1}_{\rm ang}(\mathbb R^{n})}\\
&=\sum\limits_{k=-\infty}^{\infty}2^{k\alpha q} \Bigg(\int_{2^{k-1}}^{2^k}\Big(\int_{S^{n-1}}|\mathcal{H}_{\Phi}f(r\theta)|^{p_1}d\sigma(\theta)\Big)^{\frac{p_2}{p_1}}r^{n-1}dr \Bigg)^{\frac{q}{p_2}}\\
&={\omega_{n-1}}^{\frac{1}{p_1}}\sum\limits_{k=-\infty}^{\infty}2^{k\alpha q} \Bigg(\int_{2^{k-2}}^{2^k}\Big|\sum\limits_{j=-\infty}^{\infty}\int_{E_j}\frac{\Phi(y)}{|y|^{n}}f(\frac{r}{|y|})dy\Big|^{p_2}r^{n-1}dr \Bigg)^{\frac{q}{p_2}}\\
&\leq{\omega_{n-1}}^{\frac{1}{p_1}}\sum\limits_{k=-\infty}^{\infty}2^{k\alpha q} \Bigg(\sum\limits_{j=-\infty}^{\infty}\Big(\int_{2^{k-1}}^{2^k}\Big|\int_{E_j}\frac{\Phi(y)}{|y|^{n}}f(\frac{r}{|y|})dy\Big|^{p_2}r^{n-1}dr\Big)^{\frac{1}{p_2}}\Bigg)^{q}\\
&\leq{\omega_{n-1}}^{\frac{1}{p_1}}\sum\limits_{k=-\infty}^{\infty}2^{k\alpha q} \Bigg(\sum\limits_{j=-\infty}^{\infty}\int_{E_j}\frac{|\Phi(y)|}{|y|^{n}}\Big(\int_{2^{k-1}}^{2^k}|f(\frac{r}{|y|})|^{p_2}r^{n-1}dr\Big)^{\frac{1}{p_2}}dy\Bigg)^{q}.
\end{align*}
Note that $y\in E_j$,
\begin{align*}
&\int_{E_j}\frac{|\Phi(y)|}{|y|^{n}}\Big(\int_{2^{k-1}}^{2^k}|f(\frac{r}{|y|})|^{p_2}r^{n-1}dr\Big)^{\frac{1}{p_2}}dy\\
&=\int_{E_j}\frac{|\Phi(y)|}{|y|^{n}}|y|^{\frac{n}{p_2}}\Big(\int_{2^{k-1}}^{2^k}|f(r)|^{p_2}r^{n-1}dr\Big)^{\frac{1}{p_2}}dy\\
&\backsimeq{\omega_{n-1}}^{\frac{1}{p_1}-\frac{1}{\tilde{p_1}}}\int_{E_j}\frac{|\Phi(y)|}{|y|^{n}}|y|^{\frac{n}{p_2}}\left\|f\chi_{E_{k-j}}\right\|_{L^{p_2}_{rad}L^{\tilde{p_1}}_{ang}\left(\mathbb R^n\right)}dy.
\end{align*}
Thus, we have
\begin{align*}
\left\|\mathcal{H}_{\Phi}f\right\|_{{\dot{K}}^{\alpha,q}_{L^{p_2}_{\rm rad}L^{p_1}_{\rm ang}}\left(\mathbb R^n\right)}&\lesssim\sum\limits_{j=-\infty}^{\infty}\Big(\int_{E_j}\frac{|\Phi(y)|}{|y|^{n}}|y|^{\frac{n}{p_2}+\alpha}\Big(\sum\limits_{k=-\infty}^{\infty}2^{(k-j)\alpha q} \left\|f\chi_{E_{k-j}}\right\|^q_{L^{p_2}_{rad}L^{\tilde{p_1}}_{ang}\left(\mathbb R^n\right)}\Big)^{1/q}dy\Big)\\
&\lesssim\int_{\mathbb R^{n}}\frac{|\Phi(y)|}{|y|^{n}}|y|^{\frac{n}{p_2}+\alpha}dy\left\|f\right\|_{{\dot{K}}^{\alpha,q}_{L^{p_2}_{\rm rad}L^{\tilde{p_1}}_{\rm ang}}\left(\mathbb R^n\right)}.
\end{align*}
Therefore, Theorem 1.5 is finished.
\end{proof}
\section{Mixed radial-angular integrabilities for fractional Hausdorff operator }

The $n$-dimensional fractional Hausdorff operator for a general function $\Phi$ is defined by

$$H_{\Phi,\beta}f(x)=\int_{\mathbb{R}^n}\frac{\Phi(x|y|^{-1})}{|y|^{n-\beta}}f(y)dy, \quad 0\leq\beta<n.$$

\noindent\textbf{Theorem 3.1} Let $0\leq\beta<n$, $\alpha, \gamma\in\mathbb R $, $1\leq p_1, p_2, q <\infty$, $q\leq p_2$, $\frac{\gamma+n}{p_2}=\frac{\alpha+n}{p_1}-\beta$ and $\frac{1}{p_2}+1=\frac{1}{s}+\frac{1}{p_1}$, then
$$(a) \quad \|H_{\Phi,\beta}f\|_{L^{p_2}_{rad}L^{q}_{ang}(\mathbb R^{n}, r^\gamma)}\leq \|\Phi\|_{L^{q}_{ang}L^{s}_{rad}(\mathbb R^{n}, r^{\frac{(\gamma+n)s}{p_2}-n})}\|f\|_{L^{1}_{ang}L^{p_1}_{rad}(\mathbb R^{n}, r^\alpha)};$$

$$(b) \quad\|H_{\Phi,\beta}f\|_{L^{p_2}_{rad}L^{q}_{ang}(\mathbb R^{n}, r^\gamma)}\leq \|\Phi\|_{L^{q}_{ang}L^{s}_{rad}(\mathbb R^{n}, r^{\frac{(\gamma+n)s}{p_2}-n})}\|f\|_{L^{p_1}_{rad}L^{q}_{ang}(\mathbb R^{n}, r^\alpha)}.$$

\begin{proof}[\textbf{Proof of Theorem 3.1 (a)}] By polar coordinates, we obtain that

$$H_{\Phi,\beta}f(x)=\int_{0}^{\infty}\int_{S^{n-1}}{\Phi(\frac{x}{t})}t^\beta f(t\theta_1)d\sigma(\theta_1)\frac{dt}{t}.$$
Thus, we have
\begin{align*}
&\|H_{\Phi,\beta}f\|_{L^{p_2}_{rad}L^{q}_{ang}(\mathbb R^{n}, r^\gamma)}\\
&=\Big(\int_{0}^{\infty}\Big(\int_{S^{n-1}}\Big|\int_{0}^{\infty}\int_{S^{n-1}}{\Phi(\frac{r\theta}{t})}t^\beta f(t\theta_1)d\sigma(\theta_1)\frac{dt}{t}\Big|^{q}d\sigma(\theta)\Big)^{\frac{p_2}{q}}r^{n+\gamma-1}dr\Big)^{\frac{1}{p_2}}.
\end{align*}
Applying the conditions $\frac{\gamma+n}{p_2}=\frac{\alpha+n}{p_1}-\beta$ and $q\leq p_2$ give rise to
\begin{align*}
&\Big(\int_{0}^{\infty}\Big(\int_{S^{n-1}}\Big|\int_{0}^{\infty}\int_{S^{n-1}}{\Phi(\frac{r\theta}{t})}t^\beta f(t\theta_1)d\sigma(\theta_1)\frac{dt}{t}\Big|^{q}d\sigma(\theta)\Big)^{\frac{p_2}{q}}r^{n+\gamma-1}dr\Big)^{\frac{1}{p_2}}\\
&=\Big(\int_{0}^{\infty}\Big(\int_{S^{n-1}}\Big|\int_{0}^{\infty}\int_{S^{n-1}}{\Phi(\frac{r\theta}{t})}(\frac{r}{t})^{\frac{\gamma+n}{p_2}} f(t\theta_1)t^{\frac{\alpha+n}{p_1}}d\sigma(\theta_1)\frac{dt}{t}\Big|^{q}d\sigma(\theta)\Big)^{\frac{p_2}{q}}\frac{dr}{r}\Big)^{\frac{1}{p_2}}\\
&\leq\Big(\int_{0}^{\infty}\Big(\int_{S^{n-1}}\Big(\int_{0}^{\infty}\int_{S^{n-1}}\big|{\Phi(\frac{r\theta}{t})}\big|(\frac{r}{t})^{\frac{\gamma+n}{p_2}} \big|f(t\theta_1)\big|t^{\frac{\alpha+n}{p_1}}d\sigma(\theta_1)\frac{dt}{t}\Big)^{q}d\sigma(\theta)\Big)^{\frac{p_2}{q}}\frac{dr}{r}\Big)^{\frac{1}{p_2}}\\
&=\Big(\int_{0}^{\infty}\Big(\int_{S^{n-1}}\Big(\int_{S^{n-1}}\int_{0}^{\infty}\big|{\Phi(\frac{r\theta}{t})}\big|(\frac{r}{t})^{\frac{\gamma+n}{p_2}} \big|f(t\theta_1)\big|t^{\frac{\alpha+n}{p_1}}d\sigma(\theta_1)\frac{dt}{t}\Big)^{q}d\sigma(\theta)\Big)^{\frac{p_2}{q}}\frac{dr}{r}\Big)^{\frac{1}{p_2}}\\
&\leq\Big(\int_{S^{n-1}}\Big(\int_{0}^{\infty}\Big(\int_{S^{n-1}}\int_{0}^{\infty}\big|{\Phi(\frac{r\theta}{t})}\big|(\frac{r}{t})^{\frac{\gamma+n}{p_2}} \big|f(t\theta_1)\big|t^{\frac{\alpha+n}{p_1}}d\sigma(\theta_1)\frac{dt}{t}\Big)^{p_2}\frac{dr}{r}\Big)^{\frac{q}{p_2}}d\sigma(\theta)\Big)^{\frac{1}{q}}.
\end{align*}
By the Minkowski inequality, we get
\begin{align*}
&\Big(\int_{0}^{\infty}\Big(\int_{S^{n-1}}\int_{0}^{\infty}\big|{\Phi(\frac{r\theta}{t})}\big|(\frac{r}{t})^{\frac{\gamma+n}{p_2}} \big|f(t\theta_1)\big|t^{\frac{\alpha+n}{p_1}}d\sigma(\theta_1)\frac{dt}{t}\Big)^{p_2}\frac{dr}{r}\Big)^{\frac{1}{p_2}}\\
&\leq\int_{S^{n-1}}\Big(\int_{0}^{\infty}\Big|\int_{0}^{\infty}\big|{\Phi(\frac{r\theta}{t})}\big|(\frac{r}{t})^{\frac{\gamma+n}{p_2}} \big|f(t\theta_1)\big|t^{\frac{\alpha+n}{p_1}}\frac{dt}{t}\Big|^{p_2}\frac{dr}{r}\Big)^{\frac{1}{p_2}}d\sigma(\theta_1).
\end{align*}
Note that $\frac{1}{p_2}+1=\frac{1}{s}+\frac{1}{p_1}$, by Young's inequality for convolution, we have
\begin{align*}
&\Big(\int_{0}^{\infty}\Big|\int_{0}^{\infty}\big|{\Phi(\frac{r\theta}{t})}\big|(\frac{r}{t})^{\frac{\gamma+n}{p_2}} \big|f(t\theta_1)\big|t^{\frac{\alpha+n}{p_1}}\frac{dt}{t}\Big|^{p_2}\frac{dr}{r}\Big)^{\frac{1}{p_2}}\\
&\leq\Big(\int_{0}^{\infty}\big|{\Phi({r\theta})}\big|^sr^{\frac{(\gamma+n)s}{p_2}}\frac{dr}{r}\Big)^{\frac{1}{s}}
\Big(\int_{0}^{\infty}\big|f(r\theta_1)\big|^{p_1}r^{\alpha+n}\frac{dr}{r}\Big)^{\frac{1}{p_1}}\\
&=\Big(\int_{0}^{\infty}\big|{\Phi(t\theta)}\big|^st^{\frac{(\gamma+n)s}{p_2}-1}dt\Big)^{\frac{1}{s}}
\Big(\int_{0}^{\infty}\big|f(t\theta_1)\big|^{p_1}t^{\alpha+n-1}dt\Big)^{\frac{1}{p_1}}.
\end{align*}
Therefore,
\begin{align*}
&\|H_{\Phi,\beta}f\|_{L^{p_2}_{rad}L^{q_2}_{ang}(\mathbb R^{n}, r^\gamma)}\\
&\leq\Big\{\int_{S^{n-1}}\Big[\int_{S^{n-1}}\Big(\int_{0}^{\infty}\big|{\Phi(t\theta)}\big|^st^{\frac{(\gamma+n)s}{p_2}-1}dt\Big)^{\frac{1}{s}}
\Big(\int_{0}^{\infty}\big|f(t\theta_1)\big|^{p_1}t^{\alpha+n-1}dt\Big)^{\frac{1}{p_1}}d\sigma(\theta_1)\Big]^{q}d\sigma(\theta)\Big\}^{\frac{1}{q}}\\
&=\Big(\int_{S^{n-1}}\Big(\int_{0}^{\infty}\big|{\Phi(t\theta)}\big|^st^{\frac{(\gamma+n)s}{p_2}-1}dt\Big)^{\frac{q}{s}}d\sigma(\theta)\Big)^{\frac{1}{q}}
\int_{S^{n-1}}\Big(\int_{0}^{\infty}\big|f(t\theta_1)\big|^{p_1}t^{\alpha+n-1}dt\Big)^{\frac{1}{p_1}}d\sigma(\theta_1)\\
&\leq\|\Phi\|_{L^{q}_{ang}L^{s}_{rad}(\mathbb R^{n}, r^{\frac{(\gamma+n)s}{p_2}-n})}\|f\|_{L^{1}_{ang}L^{p_1}_{rad}(\mathbb R^{n}, r^\alpha)}.
\end{align*}
Theorem 3.1 $(a)$ is completed.
\end{proof}

\begin{proof}[\textbf{Proof of Theorem 3.1 (b)}] Adopting the same method as Theorem 1.4, we can restrict to $f$ radial function.
Thus, by using the conditions $q\leq p_2$ and  $\frac{\gamma+n}{p_2}=\frac{\alpha+n}{p_1}-\beta$, we derive that
\begin{align*}
&\|H_{\Phi,\beta}f\|_{L^{p_2}_{rad}L^{q}_{ang}(\mathbb R^{n}, r^\gamma)}\\
&=\Big(\int_{0}^{\infty}\Big(\int_{S^{n-1}}\Big|\int_{0}^{\infty}\int_{S^{n-1}}{\Phi(\frac{r\theta}{t})}t^{\beta-1} f(t\theta_1)d\sigma(\theta_1)dt\Big|^{q}d\sigma(\theta)\Big)^{\frac{p_2}{q}}r^{n+\gamma-1}dr\Big)^{\frac{1}{p_2}}\\
&=\omega_{n-1}^{\frac{1}{q}}\int_{0}^{\infty}\Big(\int_{S^{n-1}}\Big|\int_{0}^{\infty}{\Phi(\frac{r\theta}{t})}t^{\beta-1} f(t)dt\Big|^{q}d\sigma(\theta)\Big)^{\frac{p_2}{q}}r^{n+\gamma-1}dr\Big)^{\frac{1}{p_2}}\\
&\leq\omega_{n-1}^{\frac{1}{q}}\Big(\int_{S^{n-1}}\Big(\int_{0}^{\infty}\Big|\int_{0}^{\infty}{\Phi(\frac{r\theta}{t})}(\frac{r}{t})^{\frac{\gamma+n}{p_2}} f(t)t^{\frac{\alpha+n}{p_1}}\frac{dt}{t}\Big|^{p_2}\frac{dr}{r}\Big)^{\frac{q}{p_2}}d\sigma(\theta)\Big)^{\frac{1}{q}}.
\end{align*}
Note that $\frac{1}{p_2}+1=\frac{1}{s}+\frac{1}{p_1}$, by Young's inequality for convolution, we have
\begin{align*}
&\Big(\int_{0}^{\infty}\Big|\int_{0}^{\infty}\big|{\Phi(\frac{r\theta}{t})}\big|(\frac{r}{t})^{\frac{\gamma+n}{p_2}} \big|f(t)\big|t^{\frac{\alpha+n}{p_1}}\frac{dt}{t}\Big|^{p_2}\frac{dr}{r}\Big)^{\frac{1}{p_2}}\\
&\leq\Big(\int_{0}^{\infty}\big|{\Phi(t\theta)}\big|^st^{\frac{(\gamma+n)s}{p_2}-1}dt\Big)^{\frac{1}{s}}
\Big(\int_{0}^{\infty}\big|f(t)\big|^{p_1}t^{\alpha+n-1}dt\Big)^{\frac{1}{p_1}}.
\end{align*}
Therefore,
\begin{align*}
&\|H_{\Phi,\beta}f\|_{L^{p_2}_{rad}L^{q_2}_{ang}(\mathbb R^{n}, r^\gamma)}\\
&\leq\omega_{n-1}^{\frac{1}{q}}\Big\{\int_{S^{n-1}}\Big[\Big(\int_{0}^{\infty}\big|{\Phi(t\theta)}\big|^st^{\frac{(\gamma+n)s}{p_2}-1}dt\Big)^{\frac{1}{s}}
\Big(\int_{0}^{\infty}\big|f(t)\big|^{p_1}t^{\alpha+n-1}dt\Big)^{\frac{1}{p_1}}\Big]^{q}d\sigma(\theta)\Big\}^{\frac{1}{q}}\\
&=\omega_{n-1}^{\frac{1}{q}}\Big(\int_{S^{n-1}}\Big(\int_{0}^{\infty}\big|{\Phi(t\theta)}\big|^st^{\frac{(\gamma+n)s}{p_2}-1}dt\Big)^{\frac{q}{s}}d\sigma(\theta)\Big)^{\frac{1}{q}}
\Big(\int_{0}^{\infty}\big|f(t\theta_1)\big|^{p_1}t^{\alpha+n-1}dt\Big)^{\frac{1}{p_1}}\\
&=\|\Phi\|_{L^{q}_{ang}L^{s}_{rad}(\mathbb R^{n}, r^{\frac{(\gamma+n)s}{p_2}-n})}\|f\|_{L^{p_1}_{rad}L^{q}_{ang}(\mathbb R^{n}, r^\alpha)}.
\end{align*}
Theorem 3.1 $(b)$ is finished.
\end{proof}

Next, we will establish mixed radial-angular weak bounds for fractional Hausdorff operator.

\noindent\textbf{Theorem 3.2} Suppose that $0\leq\beta<n$, $1\leq p_1, p_2, q<\infty$, $\alpha\in\mathbb R $, $n+\gamma>0 $ and $\frac{\gamma+n}{p_2}=\frac{\alpha+n}{p_1}-\beta$, then

$$\|H_{\Phi,\beta}f\|_{\mathcal{W}L^{p_2}_{rad}L^{q}_{ang}(\mathbb R^{n},r^\gamma)}\leq\frac{{{\omega}_{n}}^{1/q}}{(n+\gamma)^{1/p_2}}\|\Phi\|_ {L^{\infty}_{ang}L^{p_1}_{rad}(\mathbb R^{n},r^{\frac{(\gamma+n)p'_1}{p_2}-n})}\|f\|_ {L^{p_1}_{rad}L^{1}_{ang}(\mathbb R^{n},r^\alpha)}.$$

\begin{proof} By using spherical coordinates, we have
\begin{align*}
|H_{\Phi,\beta}f(x)|&=\Big|\int_{\mathbb{R}^n}\frac{\Phi(x|y|^{-1})}{|y|^{n-\beta}}f(y)dy\Big|\\
&=\Big|\int_{\mathbb{R}^n}\frac{\Phi(x|y|^{-1})}{|y|^{n-\beta}}f(y)|y|^{\frac{\alpha}{p_1}}|y|^{-\frac{\alpha}{p_1}}dy\Big|\\
&=\Big|\int_{0}^{\infty}\int_{S^{n-1}}\Phi(xt^{-1})t^{\beta-n-\frac{\alpha}{p_1}}f(t\theta_1)t^{\frac{\alpha}{p_1}}d\sigma(\theta_1)t^{n-1}dt\Big|\\
&=\int_{0}^{\infty}|\Phi(xt^{-1})|\Big(\int_{S^{n-1}}|f(t\theta_1)|d\sigma(\theta_1)\Big)t^{\beta-\frac{\alpha}{p_1}-1}t^{\frac{\alpha}{p_1}}d\sigma(\theta_1)dt\\
&\leq\int_{0}^{\infty}|\Phi(xt^{-1})|t^{(\beta-n-\frac{\alpha}{p_1})p'_1+n-1}dt\Big)^{\frac{1}{p'_1}}\|f\|_ {L^{p_1}_{rad}L^{1}_{ang}(\mathbb R^{n},r^\alpha)}\\
&=|x|^{-\frac{\gamma+n}{p_2}}\int_{0}^{\infty}|\Phi(x't^{-1})|t^{\frac{(\gamma+n)p'_1}{p_2}-1}dt\Big)^{\frac{1}{p'_1}}\|f\|_ {L^{p_1}_{rad}L^{1}_{ang}(\mathbb R^{n},r^\alpha)}\\
&\leq|x|^{-\frac{\gamma+n}{p_2}}\|\Phi\|_ {L^{\infty}_{ang}L^{p_1}_{rad}(\mathbb R^{n},r^{\frac{(\gamma+n)p'_1}{p_2}-n})}\|f\|_ {L^{p_1}_{rad}L^{1}_{ang}(\mathbb R^{n},r^\alpha)}\\
&=K|x|^{-\frac{\gamma+n}{p_2}}.
\end{align*}
Thus,
\begin{align*}
&\big\|\chi_{\{x\in\mathbb R^n:  |H_{\Phi,\beta}f(x)|>\lambda\}}\big\|_{L^{p_2}_{rad}L^{q}_{ang}(\mathbb R^n,r^\gamma)}\\
&\leq\Big\|\chi_{\{x\in\mathbb R^n:  K|x|^{-\frac{\gamma+n}{p_2}}>\lambda\}}\Big\|_{L^{p_2}_{rad}L^{q}_{ang}(\mathbb R^n,r^\gamma)}\\
&=\Big\|\chi_{\{x\in\mathbb R^n: |x|<\big(\frac{K}{\lambda}\big)^{\frac{p_2}{\gamma+n}}\}}\Big\|_{L^{p_2}_{rad}L^{q}_{ang}(\mathbb R^n,r^\gamma)}\\
&={{\omega}_{n}}^{1/q}\Bigg(\int_{0}^{\big(\frac{K}{\lambda}\big)^{\frac{p_2}{\gamma+n}}}r^{n+\gamma-1}dr\Bigg)^{1/p_2}\\
&=\frac{{{\omega}_{n}}^{1/q}}{(n+\gamma)^{1/p_2}\lambda}\|\Phi\|_ {L^{\infty}_{ang}L^{p_1}_{rad}(\mathbb R^{n},r^{\frac{(\gamma+n)p'_1}{p_2}-n})}\|f\|_ {L^{p_1}_{rad}L^{1}_{ang}(\mathbb R^{n},r^\alpha)}.
\end{align*}
This implies Theorem 3.2.
\end{proof}
\section{Mixed radial-angular integrabilities for multilinear  Hausdorff operator }

The multilinear  Hausdorff operator is defined by
$$\mathcal{S}_{\Psi}(f_1,f_2,\ldots,f_m)(x)=\int_{\mathbb R^{nm}}\frac{\Psi(\frac{x}{|y_1|},\frac{x}{|y_2|},\ldots,\frac{x}{|y_m|})}{|y_1|^n|y_2|^n\ldots|y_m|^n}\prod_{i=1}^{m}f_i(y_i)dy_1y_2\ldots dy_m,$$
where $\Psi(s_1,s_2,\ldots,s_m)$ is a locally integrable radial function on $\mathbb R^{+}\times\mathbb R^{+}\times\ldots\times\mathbb R^{+}$.

\noindent\textbf{Theorem 4.1} Let $\alpha=\sum_{i=1}^{m}\alpha_i$, $1\geq\frac{1}{p}=\sum_{i=1}^{m}\frac{1}{p_i}$, $1\geq\frac{1}{q}=\sum_{i=1}^{m}\frac{1}{q_i}$, $1\leq\tilde{p},\tilde{p}_i<\infty$. If
$$K_0=\omega_{n-1}^{m+\frac{1}{\tilde{p}}-\sum_{i=1}^{m}\frac{1}{\tilde{p}}_i}\int_{0}^{\infty}\int_{0}^{\infty}\cdots\int_{0}^{\infty}\frac{\Psi(r_1,r_2,\ldots,r_m)}{\prod_{i=1}^{m}r_i^{1-n/p_i-\alpha_i}}dr_1dr_2\ldots dr_m<\infty.$$
Then
$$
\left\|\mathcal{S}_{\Psi}(f_1,f_2,\ldots,f_m)\right\|_{{\dot{K}}^{\alpha,q}_{L^{p}_{\rm rad}L^{\tilde{p}}_{\rm ang}}\left(\mathbb R^n\right)}\leq K_0\prod_{i=1}^{m}\left\|f_i\right\|_{{\dot{K}}^{\alpha_i,q_i}_{L^{p_i}_{\rm rad}L^{\tilde{p}_i}_{\rm ang}}\left(\mathbb R^n\right)}.
$$

\begin{proof}By polar coordinates, we obtain that
\begin{align*}
\mathcal{S}_{\Psi}(f_1,f_2)(x)&=\int_{0}^{\infty}\int_{0}^{\infty}\frac{\Psi(\frac{|x|}{t_1},\frac{|x|}{t_2})}{t_1t_2}\prod_{i=1}^{2}\Big(\int_{S^{n-1}}f_i(t_iy'_i)d\sigma(y'_i)\Big)dt_1dt_2\\
&=\int_{0}^{\infty}\int_{0}^{\infty}\frac{\Psi(r_1,r_2)}{r_1r_2}\prod_{i=1}^{2}\Big(\int_{S^{n-1}}f_i(|x|r_{i}^{-1}y'_i)d\sigma(y'_i)\Big)dr_1dr_2.
\end{align*}
Applying the Minkowski inequality, we have
\begin{align*}
&\|\mathcal{S}_{\Psi}(f_1,f_2)\chi_k\|_{L^{p}_{\rm rad}L^{\tilde{p}}_{\rm ang}(\mathbb R^{n})}\\
&\leq\omega_{n-1}^{\frac{1}{\tilde{p}}}\int_{0}^{\infty}\int_{0}^{\infty}\frac{\Psi(r_1,r_2)}{r_1r_2}\Big\|\prod_{i=1}^{2}\Big(\int_{S^{n-1}}f_i(rr_{i}^{-1}y'_i)d\sigma(y'_i)\Big)\chi_k\Big\|_{L^{p}(\mathbb R^{+},r^{n-1})}dr_1dr_2\\
&\leq\omega_{n-1}^{\frac{1}{\tilde{p}}}\int_{0}^{\infty}\int_{0}^{\infty}\frac{\Psi(r_1,r_2)}{r_1r_2}\prod_{i=1}^{2}\Big\|\Big(\int_{S^{n-1}}f_i(rr_{i}^{-1}y'_i)d\sigma(y'_i)\Big)\chi_k\Big\|_{L^{p_i}(\mathbb R^{+},r^{n-1})}dr_1dr_2\\
&=\omega_{n-1}^{\frac{1}{\tilde{p}}-\sum_{i=1}^{2}\frac{1}{\tilde{p}}_i}\int_{0}^{\infty}\int_{0}^{\infty}\frac{\Psi(r_1,r_2)}{r_1r_2}\prod_{i=1}^{2}\Big\|\Big(\int_{S^{n-1}}f_i(|x|r_{i}^{-1}y'_i)d\sigma(y'_i)\Big)\chi_k\Big\|_{L^{p_i}_{\rm rad}L^{\tilde{p}_i}_{\rm ang}(\mathbb R^{n})}dr_1dr_2.
\end{align*}
Note that
\begin{align*}
&\Big\|\Big(\int_{S^{n-1}}f_i(|x|r_{i}^{-1}y'_i)d\sigma(y'_i)\Big)\chi_k\Big\|_{L^{p_i}_{\rm rad}L^{\tilde{p}_i}_{\rm ang}(\mathbb R^{n})}\\
&=r_{i}^{n/p_i}\Big\|\Big(\int_{S^{n-1}}f_i(|x|y'_i)d\sigma(y'_i)\Big)\chi_{k/r_{i}}\Big\|_{L^{p_i}_{\rm rad}L^{\tilde{p}_i}_{\rm ang}(\mathbb R^{n})}\\
&\leq r_{i}^{n/p_i}\omega_{n-1}^{\frac{1}{p'_i}}\Big\|\Big(\int_{S^{n-1}}\Big|f_i(|x|y'_i)\Big|^{p_i}d\sigma(y'_i)\Big)^{\frac{1}{p_i}}\chi_{k/r_{i}}\Big\|_{L^{p_i}_{\rm rad}L^{\tilde{p}_i}_{\rm ang}(\mathbb R^{n})}\\
&=r_{i}^{n/p_i}\omega_{n-1}^{\frac{1}{p'_i}+\frac{1}{p_i}+\sum_{i=1}^{2}\frac{1}{\tilde{p}_i}}\Big(\int_{0}^{\infty}\int_{S^{n-1}}\Big|f_i(ry'_i)\Big|^{p_i}d\sigma(y'_i)\chi_{2^{k-1}<r_ir\leq2^k(r)}r^{n-1}dr\Big)^{\frac{1}{p_i}}\\
&=r_{i}^{n/p_i}\omega_{n-1}\|f_i\chi_{k/r_i}\|_{L^{p_i}_{\rm rad}L^{\tilde{p}_i}_{\rm ang}(\mathbb R^{n})}.
\end{align*}
This implies that
\begin{align*}
&\|\mathcal{S}_{\Psi}(f_1,f_2)\chi_k\|_{L^{p}_{\rm rad}L^{\tilde{p}}_{\rm ang}(\mathbb R^{n})}\\
&=\omega_{n-1}^{2+\frac{1}{\tilde{p}}-\sum_{i=1}^{2}\frac{1}{\tilde{p}}_i}\int_{0}^{\infty}\int_{0}^{\infty}\frac{\Psi(r_1,r_2)}{r_1^{1-n/p_1}r_2^{1-n/p_2}}\prod_{i=1}^{2}\|f_i\chi_{k/r_i}\|_{L^{p_i}_{\rm rad}L^{\tilde{p}_i}_{\rm ang}(\mathbb R^{n})}dr_1dr_2.
\end{align*}
Thus,
\begin{align*}
&\left\|\mathcal{S}_{\Psi}(f_1,f_2)\right\|_{{\dot{K}}^{\alpha,q}_{L^{p}_{\rm rad}L^{\tilde{p}}_{\rm ang}}\left(\mathbb R^n\right)}\\
&=\omega_{n-1}^{2+\frac{1}{\tilde{p}}-\sum_{i=1}^{2}\frac{1}{\tilde{p}}_i}\int_{0}^{\infty}\int_{0}^{\infty}\frac{\Psi(r_1,r_2)}{r_1^{1-n/p_1}r_2^{1-n/p_2}}\prod_{i=1}^{2}\Big(\sum\limits_{k=-\infty}^{\infty}2^{k\alpha_i q_i} \|f_i\chi_{k/r_i}\|^{q_i}_{L^{p_2}_{\rm rad}L^{p_1}_{\rm ang}(\mathbb R^{n})}\Big)^{\frac{1}{q_i}}dr_1dr_2\\
&=\omega_{n-1}^{2+\frac{1}{\tilde{p}}-\sum_{i=1}^{2}\frac{1}{\tilde{p}}_i}\int_{0}^{\infty}\int_{0}^{\infty}\frac{\Psi(r_1,r_2)}{r_1^{1-n/p_1-\alpha_1}r_2^{1-n/p_2-\alpha_2}}dr_1dr_2\prod_{i=1}^{2}\left\|f_i\right\|_{{\dot{K}}^{\alpha_i,q_i}_{L^{p_i}_{\rm rad}L^{\tilde{p}_i}_{\rm ang}}\left(\mathbb R^n\right)}.
\end{align*}
Theorem 4.1 is finished.
\end{proof}

\noindent\textbf{Theorem 4.2} Let $1\leq p, p_i, \tilde{p}<\infty$, $\lambda=\sum_{i=1}^{m}\lambda_i$, $\frac{1}{p}=\sum_{i=1}^{m}\frac{1}{p_i}$. If $f$ is a radial function and
$$K_1=\omega_{n-1}^{\frac{1}{\tilde{p}}-\sum_{i=1}^{m}\frac{1}{\tilde{p}}_i}\int_{0}^{\infty}\int_{0}^{\infty}\cdots\int_{0}^{\infty}\frac{\Psi(r_1,r_2,\ldots,r_m)}{\prod_{i=1}^{m}r_i^{1+n\lambda_i}}dr_1dr_2\ldots dr_m<\infty.$$
Then
$$
\left\|\mathcal{S}_{\Psi}(f_1,f_2,\ldots,f_m)\right\|_{\mathcal{B}{L^{p,\lambda}_{\rm rad}L^{\tilde{p}}_{\rm ang}}(\mathbb R^{n})}\leq K_1\prod_{i=1}^{m}\left\|f_i\right\|_{\mathcal{B}{L^{p_i,\lambda_i}_{\rm rad}L^{\tilde{p}_i}_{\rm ang}}(\mathbb R^{n})}.
$$

\begin{proof}By the Minkowski inequality and  H\"{o}lder's inequality, we have
\begin{align*}
&\Big(\frac{1}{\nu_nr^{n+n\lambda p}}\int_{0}^{r}\Big(\int_{{S}^{n-1}}\Big|\int_{0}^{\infty}\int_{0}^{\infty}\frac{\Psi(r_1,r_2)}{r_1r_2}\\
&\quad\times\prod_{i=1}^{2}\Big(\int_{S^{n-1}}f_i(|x|r_{i}^{-1}y'_i)d\sigma(y'_i)\Big)dr_1dr_2\Big|^{\tilde{p}}d\sigma(\theta)\Big)^{p/\tilde{p}}\rho^{n-1}d\rho\Big)^{1/p}\\
&=\omega_{n-1}^{\frac{1}{\tilde{p}}}\Big(\frac{1}{\nu_nr^{n+n\lambda p}}\int_{0}^{r}\Big|\int_{0}^{\infty}\int_{0}^{\infty}\frac{\Psi(r_1,r_2)}{r_1r_2}\prod_{i=1}^{2}\Big(\int_{S^{n-1}}f_i(rr_{i}^{-1}y'_i)d\sigma(y'_i)\Big)dr_1dr_2\Big|^p\rho^{n-1}d\rho\Big)^{1/p}\\
&\leq\omega_{n-1}^{\frac{1}{\tilde{p}}}\int_{0}^{\infty}\int_{0}^{\infty}\frac{\Psi(r_1,r_2)}{r_1r_2}\Big(\frac{1}{\nu_nr^{n+n\lambda p}}\int_{0}^{r}\Big|\prod_{i=1}^{2}\int_{S^{n-1}}f_i(\rho r_{i}^{-1}y'_i)d\sigma(y'_i)\Big|^p\rho^{n-1}d\rho\Big)^{1/p}dr_1dr_2\\
&\leq\omega_{n-1}^{\frac{1}{\tilde{p}}}\int_{0}^{\infty}\int_{0}^{\infty}\frac{\Psi(r_1,r_2)}{r_1r_2}\prod_{i=1}^{2}\Big(\frac{1}{\nu_nr^{n+n\lambda_i p_i}}\int_{0}^{r}\Big|\int_{S^{n-1}}f_i(\rho r_{i}^{-1}y'_i)d\sigma(y'_i)\Big|^{p_i}\rho^{n-1}d\rho\Big)^{1/p_i}dr_1dr_2\\
&=\omega_{n-1}^{\frac{1}{\tilde{p}}}\int_{0}^{\infty}\int_{0}^{\infty}\frac{\Psi(r_1,r_2)}{r_1^{1+n\lambda_1}r_2^{1+n\lambda_2}}\prod_{i=1}^{2}\Big(\frac{1}{\nu_n(r/r_i)^{n+n\lambda_i p_i}}\int_{0}^{r/r_i}\Big|\int_{S^{n-1}}f_i(\rho y'_i)d\sigma(y'_i)\Big|^{p_i}\rho^{n-1}d\rho\Big)^{1/p_i}dr_1dr_2\\
&=\omega_{n-1}^{\frac{1}{\tilde{p}}-\sum_{i=1}^{2}\frac{1}{\tilde{p}}_i}\int_{0}^{\infty}\int_{0}^{\infty}\frac{\Psi(r_1,r_2)}{r_1^{1+n\lambda_1}r_2^{1+n\lambda_2}}\prod_{i=1}^{2}\left\|f_i\right\|_{\mathcal{B}{L^{p_i,\lambda_i}_{\rm rad}L^{\tilde{p}_i}_{\rm ang}}(\mathbb R^{n})}.
\end{align*}
Theorem 4.1 is proved.
\end{proof}
\section{Boundedness of commutators of  Hausdorff operators on generalized  mixed radial-angular  central Morrey spaces}
In this section, our main focus is to study the boundedness of commutators of  Hausdorff operators on generalized  mixed radial-angular Morrey spaces. To achieve this, we first give the definition of the generalized  mixed radial-angular Morrey spaces as follows:

\noindent\textbf{Definition 5.1}
Let $1<p, \tilde{p}<\infty$. A function $f\in L^{p}_{\rm rad,{loc}}L^{\tilde{p}}_{\rm ang}((0,\infty)$ $\times{S}^{n-1})$ is said to belong to the  generalized  mixed radial-angular  central Morrey spaces ${\mathcal{\varphi}L^{p_2}_{\rm rad}L^{p_1}_{\rm ang}(\mathbb R^n)}$, if
\begin{small}
$$\|f\|_{{\mathcal{\varphi}L^{p}_{\rm rad}L^{\tilde{p}}_{\rm ang}(\mathbb R^n)}}=\sup\limits_{r>0}\frac{1}{\varphi(r)}\Big(\frac{1}{\nu_nr^{n}}\int_{0}^{r}\Big(\int_{{S}^{n-1}}|f(\rho\theta)|^{\tilde{p}}d\sigma(\theta)\Big)^{p/\tilde{p}}\rho^{n-1}d\rho\Big)^{1/p}<\infty.$$
\end{small}

Next, we will consider the boundedness  of  Hausdorff operators $\mathcal{H}_{\Phi}$
on  generalized  mixed radial-angular Morrey spaces  ${\mathcal{\varphi}L^{p}_{\rm rad}L^{\tilde{p}}_{\rm ang}(\mathbb R^n)}$.

\noindent\textbf{Theorem 5.1}  Let $1\leq p, \tilde{p}<\infty$, if
$$\mathcal{K}=\sup\limits_{r>0}\frac{1}{\varphi(r)}\int_{\mathbb{R}^n}\frac{|\Phi(y)|}{|y|^{n}}|y|^{n/p}\varphi(r/|y|)dy<\infty.$$
Then
$$\|\mathcal{H}_{\Phi}f\|_{{\mathcal{\varphi}L^{p}_{\rm rad}L^{\tilde{p}}_{\rm ang}(\mathbb R^n)}}\leq\mathcal{K}\|f\|_{{\mathcal{\varphi}L^{p}_{\rm rad}L^{\tilde{p}}_{\rm ang}(\mathbb R^n)}}.$$
\begin{proof} The Minkowski inequality gives that
\begin{align*}
&\|\mathcal{H}_{\Phi}f\|_{{\mathcal{\varphi}L^{p}_{\rm rad}L^{\tilde{p}}_{\rm ang}(\mathbb R^n)}}\\
&=\sup\limits_{r>0}\frac{1}{\varphi(r)}\Big(\frac{1}{\nu_nr^{n}}\int_{0}^{r}\Big(\int_{{S}^{n-1}}|\int_{\mathbb{R}^n}\frac{|\Phi(y)|}{|y|^{n}}f(\frac{\rho\theta}{|y|})dy|^{\tilde{p}}d\sigma(\theta)\Big)^{p/\tilde{p}}\rho^{n-1}d\rho\Big)^{1/p}\\
&\leq\sup\limits_{r>0}\frac{1}{\varphi(r)}\Big(\frac{1}{\nu_nr^{n}}\int_{0}^{r}\Big(\int_{\mathbb{R}^n}\frac{|\Phi(y)|}{|y|^{n}}\Big(\int_{{S}^{n-1}}|f(\frac{\rho\theta}{|y|})|^{\tilde{p}}d\sigma(\theta)\Big)^{1/\tilde{p}}dy\Big)^p\rho^{n-1}d\rho\Big)^{1/p}\\
&\leq\sup\limits_{r>0}\frac{1}{\varphi(r)}\int_{\mathbb{R}^n}\frac{|\Phi(y)|}{|y|^{n}}|y|^{n/p}\Big(\frac{1}{\nu_nr^{n}}\int_{0}^{r/|y|}\Big(\int_{{S}^{n-1}}|f(\rho\theta)|^{\tilde{p}}d\sigma(\theta)\Big)^{p/\tilde{p}}\rho^{n-1}d\rho\Big)^{1/p}dy\\
&\leq\sup\limits_{r>0}\frac{1}{\varphi(r)}\int_{\mathbb{R}^n}\frac{|\Phi(y)|}{|y|^{n}}|y|^{n/p}\varphi(r/|y|)dy\|f\|_{{\mathcal{\varphi}L^{p}_{\rm rad}L^{\tilde{p}}_{\rm ang}(\mathbb R^n)}}.
\end{align*}
This finishes the proof of Theorem 5.1.
\end{proof}
Now we in this position to consider the following two types of commutators of  Hausdorff operators
$$\mathcal{H}_{\Phi}^{b}f(x)=b(x)\mathcal{H}_{\Phi}f(x)-\mathcal{H}_{\Phi}(bf)(x)=\int_{\mathbb{R}^n}\frac{\Phi(y)}{|y|^{n}}f(\frac{x}{|y|})(b(x)-b(\frac{x}{|y|}))dy$$
and
$$H_{\Phi}^bf(x)=b(x)H_{\Phi}f(x)-H_{\Phi}(bf)(x)=\int_{\mathbb{R}^n}\frac{\Phi(x|y|^{-1})}{|y|^{n}}f(y)(b(x)-b(y))dy.$$

\noindent\textbf{Theorem 5.2}  Let $1\leq p, p_i, \tilde{p}<\infty$,  $\frac{1}{p}=\frac{1}{p_1}+\frac{1}{p_2}$, $\frac{1}{\tilde{p}}=\frac{1}{\tilde{p}_1}+\frac{1}{\tilde{p}_2}$, if
\begin{align*}
\mathcal{K}_0&=\sup\limits_{r>0}\frac{1}{\varphi(r)}\int_{|y|\leq1}\frac{|\Phi(y)|}{|y|^{n}}|y|^{n/p}\varphi(r/|y|)\log\frac{2}{|y|}dy\\
&\quad+\sup\limits_{r>0}\frac{1}{\varphi(r)}\int_{|y|\geq1}\frac{|\Phi(y)|}{|y|^{n}}|y|^{n/p}\log(2|y|)\varphi(r/|y|)dy<\infty.
\end{align*}
Then
$$
\|\mathcal{H}^b_{\Phi}f\|_{{\mathcal{\varphi}L^{p}_{\rm rad}L^{\tilde{p}}_{\rm ang}(\mathbb R^n)}}\lesssim \mathcal{K}_0\|b\|_{{CMOL^{p_2}_{\rm rad}L^{\tilde{p}_2}_{\rm ang}(\mathbb R^n)}}\|f\|_{{\mathcal{\varphi}L^{p_1}_{\rm rad}L^{\tilde{p}_1}_{\rm ang}(\mathbb R^n)}}.
$$
\begin{proof} Note that
\begin{align*}
&\|\mathcal{H}^b_{\Phi}f\|_{{\mathcal{\varphi}L^{p}_{\rm rad}L^{\tilde{p}}_{\rm ang}(\mathbb R^n)}}\\
&=\sup\limits_{r>0}\frac{1}{\varphi(r)}\Big(\frac{1}{\nu_nr^{n}}\int_{0}^{r}\Big(\int_{{S}^{n-1}}|\int_{\mathbb{R}^n}\frac{\Phi(|y|)}{|y|^{n}}f(\frac{\rho\theta}{|y|})(b(\rho\theta)-b(\frac{\rho\theta}{|y|}))dy|^{\tilde{p}}d\sigma(\theta)\Big)^{p/\tilde{p}}\rho^{n-1}d\rho\Big)^{1/p}\\
&\leq\sup\limits_{r>0}\frac{1}{\varphi(r)}\Big(\frac{1}{\nu_nr^{n}}\int_{0}^{r}\Big(\int_{{S}^{n-1}}|\int_{\mathbb{R}^n}\frac{\Phi(|y|)}{|y|^{n}}f(\frac{\rho\theta}{|y|})(b(\rho\theta)-b_B)dy|^{\tilde{p}}d\sigma(\theta)\Big)^{p/\tilde{p}}\rho^{n-1}d\rho\Big)^{1/p}\\
&\quad+\sup\limits_{r>0}\frac{1}{\varphi(r)}\Big(\frac{1}{\nu_nr^{n}}\int_{0}^{r}\Big(\int_{{S}^{n-1}}|\int_{\mathbb{R}^n}\frac{\Phi(|y|)}{|y|^{n}}f(\frac{\rho\theta}{|y|})(b_B-b_{|y|^{-1}B})dy|^{\tilde{p}}d\sigma(\theta)\Big)^{p/\tilde{p}}\rho^{n-1}d\rho\Big)^{1/p}\\
&\quad+\sup\limits_{r>0}\frac{1}{\varphi(r)}\Big(\frac{1}{\nu_nr^{n}}\int_{0}^{r}\Big(\int_{{S}^{n-1}}|\int_{\mathbb{R}^n}\frac{\Phi(|y|)}{|y|^{n}}f(\frac{\rho\theta}{|y|})(b_{|y|^{-1}B}-b(\frac{\rho\theta}{|y|}))dy|^{\tilde{p}}d\sigma(\theta)\Big)^{p/\tilde{p}}\rho^{n-1}d\rho\Big)^{1/p}\\
&:=I+II+III.
\end{align*}
For $I_1$, the Minkowski inequality and  H\"{o}lder's inequality deduce that
\begin{align*}
&I\leq\sup\limits_{r>0}\frac{1}{\varphi(r)}\int_{\mathbb{R}^n}\frac{\Phi(|y|)}{|y|^{n}}\Big(\frac{1}{\nu_nr^{n}}\int_{0}^{r}\Big(\int_{{S}^{n-1}}|f(\frac{\rho\theta}{|y|})(b(\rho\theta)-b_B))|^{\tilde{p}}d\sigma(\theta)\Big)^{p/\tilde{p}}\rho^{n-1}d\rho\Big)^{1/p}dy\\
&\leq\sup\limits_{r>0}\frac{1}{\varphi(r)}\int_{\mathbb{R}^n}\frac{\Phi(|y|)}{|y|^{n}}\Big(\frac{1}{\nu_nr^{n}}\int_{0}^{r}\Big(\int_{{S}^{n-1}}|f(\frac{\rho\theta}{|y|})|^{\tilde{p}_1}d\sigma(\theta)\Big)^{1/\tilde{p}_1}\\
&\quad\times\Big(\int_{{S}^{n-1}}|b(\rho\theta)-b_B)|^{\tilde{p}_2}d\sigma(\theta)\Big)^{1/\tilde{p}_2}\Big)^p\rho^{n-1}d\rho\Big)^{1/p}dy\\
&\leq\sup\limits_{r>0}\frac{1}{\varphi(r)}\int_{\mathbb{R}^n}\frac{\Phi(|y|)}{|y|^{n}}\Big(\frac{1}{\nu_nr^{n}}\int_{0}^{r}\Big(\int_{{S}^{n-1}}|f(\frac{\rho\theta}{|y|})|^{\tilde{p}_1}d\sigma(\theta)\Big)^{p_1/\tilde{p}_1}\rho^{n-1}d\rho\Big)^{1/p_1}\\
&\quad\times\Big(\frac{1}{\nu_nr^{n}}\int_{0}^{r}\Big(\int_{{S}^{n-1}}|b(\rho\theta)-b_B)|^{\tilde{p}_2}d\sigma(\theta)\Big)^{p_2/\tilde{p}_2}\rho^{n-1}d\rho\Big)^{1/p_2}dy\\
&\leq\sup\limits_{r>0}\frac{1}{\varphi(r)}\int_{\mathbb{R}^n}\frac{\Phi(|y|)}{|y|^{n}}|y|^{n/p}\varphi(r/|y|)dy\|f\|_{{\mathcal{\varphi}L^{p_1}_{\rm rad}L^{\tilde{p}_1}_{\rm ang}(\mathbb R^n)}}\|b\|_{{CMOL^{p_2}_{\rm rad}L^{\tilde{p}_2}_{\rm ang}(\mathbb R^n)}}.
\end{align*}
For $II$, we divide it into two parts as follows:
\begin{align*}
II&\leq\sup\limits_{r>0}\frac{1}{\varphi(r)}\int_{\mathbb{R}^n}\frac{\Phi(|y|)}{|y|^{n}}\Big(\frac{1}{\nu_nr^{n}}\int_{0}^{r}\Big(\int_{{S}^{n-1}}|f(\frac{\rho\theta}{|y|})(b_B-b_{|y|^{-1}B})|^{\tilde{p}}d\sigma(\theta)\Big)^{p/\tilde{p}}\rho^{n-1}d\rho\Big)^{1/p}dy\\
&\leq\sup\limits_{r>0}\frac{1}{\varphi(r)}\int_{|y|\leq1}\frac{\Phi(|y|)}{|y|^{n}}\Big(\frac{1}{\nu_nr^{n}}\int_{0}^{r}\Big(\int_{{S}^{n-1}}|f(\frac{\rho\theta}{|y|})(b_B-b_{|y|^{-1}B})|^{\tilde{p}}d\sigma(\theta)\Big)^{p/\tilde{p}}\rho^{n-1}d\rho\Big)^{1/p}dy\\
&\quad+\sup\limits_{r>0}\frac{1}{\varphi(r)}\int_{|y|\geq1}\frac{\Phi(|y|)}{|y|^{n}}\Big(\frac{1}{\nu_nr^{n}}\int_{0}^{r}\Big(\int_{{S}^{n-1}}|f(\frac{\rho\theta}{|y|})(b_B-b_{|y|^{-1}B})|^{\tilde{p}}d\sigma(\theta)\Big)^{p/\tilde{p}}\rho^{n-1}d\rho\Big)^{1/p}dy\\
&:=II_1+II_2.
\end{align*}
Next, we estimate $II_1$, for any $|y|\leq1$, there exists $k\in\mathbb Z^+$ such that $2^{-k}\leq|y|<2^{-k+1}$.  The H\"{o}lders inequality gives that
\begin{align*}
II_1&\leq\sup\limits_{r>0}\frac{1}{\varphi(r)}\int_{|y|\leq1}\frac{\Phi(|y|)}{|y|^{n}}\varphi(r/|y|)|b_B-b_{|y|^{-1}B}|dy\|f\|_{{\mathcal{\varphi}L^{p_1}_{\rm rad}L^{\tilde{p}_1}_{\rm ang}(\mathbb R^n)}}\\
&=\sup\limits_{r>0}\frac{1}{\varphi(r)}\int_{|y|\leq1}\frac{\Phi(|y|)}{|y|^{n}}\varphi(r/|y|)\Big(\sum\limits_{i=0}^{k}|b_{2^{-i}B}-b_{2^{-i-1}B}|+|b_{2^{-k}B}-b_{|y|^{-1}B}|\Big)dy\|f\|_{{\mathcal{\varphi}L^{p_1}_{\rm rad}L^{\tilde{p}_1}_{\rm ang}(\mathbb R^n)}}\\
&\lesssim\sup\limits_{r>0}\frac{1}{\varphi(r)}\sum\limits_{k=0}^{\infty}\int_{2^{-k-1}\leq|y|\leq2^{-k}}\frac{\Phi(|y|)}{|y|^{n}}\varphi(r/|y|)(k+2)dy\|b\|_{{CMOL^{p_2}_{\rm rad}L^{\tilde{p}_2}_{\rm ang}(\mathbb R^n)}}\|f\|_{{\mathcal{\varphi}L^{p_1}_{\rm rad}L^{\tilde{p}_1}_{\rm ang}(\mathbb R^n)}}\\
&\leq\sup\limits_{r>0}\frac{1}{\varphi(r)}\sum\limits_{k=0}^{\infty}\int_{2^{-k-1}\leq|y|\leq2^{-k}}\frac{\Phi(|y|)}{|y|^{n}}\varphi(r/|y|)(\log2^k+2)dy\|b\|_{{CMOL^{p_2}_{\rm rad}L^{\tilde{p}_2}_{\rm ang}(\mathbb R^n)}}\|f\|_{{\mathcal{\varphi}L^{p_1}_{\rm rad}L^{\tilde{p}_1}_{\rm ang}(\mathbb R^n)}}\\
&\leq\sup\limits_{r>0}\frac{1}{\varphi(r)}\int_{|y|\leq1}\frac{\Phi(|y|)}{|y|^{n}}\varphi(r/|y|)(\log\frac{1}{|y|}+2)dy\|b\|_{{CMOL^{p_2}_{\rm rad}L^{\tilde{p}_2}_{\rm ang}(\mathbb R^n)}}\|f\|_{{\mathcal{\varphi}L^{p_1}_{\rm rad}L^{\tilde{p}_1}_{\rm ang}(\mathbb R^n)}}.
\end{align*}
Similarly, we also  have
$$II_2\lesssim\sup\limits_{r>0}\frac{1}{\varphi(r)}\int_{|y|\geq1}\frac{\Phi(|y|)}{|y|^{n}}\varphi(r/|y|)(\log|y|+2)dy\|b\|_{{CMOL^{p_2}_{\rm rad}L^{\tilde{p}_2}_{\rm ang}(\mathbb R^n)}}\|f\|_{{\mathcal{\varphi}L^{p_1}_{\rm rad}L^{\tilde{p}_1}_{\rm ang}(\mathbb R^n)}}.$$
Finally, we estimates $III$, by the Minkowski inequality and  H\"{o}lder's inequality, we obtain
\begin{align*}
&I\leq\sup\limits_{r>0}\frac{1}{\varphi(r)}\int_{\mathbb{R}^n}\frac{\Phi(|y|)}{|y|^{n}}\Big(\frac{1}{\nu_nr^{n}}\int_{0}^{r}\Big(\int_{{S}^{n-1}}|f(\frac{\rho\theta}{|y|})(b_{|y|^{-1}B}-b(\frac{\rho\theta}{|y|})|^{\tilde{p}}d\sigma(\theta)\Big)^{p/\tilde{p}}\rho^{n-1}d\rho\Big)^{1/p}dy\\
&\leq\sup\limits_{r>0}\frac{1}{\varphi(r)}\int_{\mathbb{R}^n}\frac{\Phi(|y|)}{|y|^{n}}\Big(\frac{1}{\nu_nr^{n}}\int_{0}^{r}\Big(\int_{{S}^{n-1}}|f(\frac{\rho\theta}{|y|})|^{\tilde{p}_1}d\sigma(\theta)\Big)^{1/\tilde{p}_1}\\
&\quad\times\Big(\int_{{S}^{n-1}}|b_{|y|^{-1}B}-b(\frac{\rho\theta}{|y|})|^{\tilde{p}_2}d\sigma(\theta)\Big)^{1/\tilde{p}_2}\Big)^p\rho^{n-1}d\rho\Big)^{1/p}dy\\
&\leq\sup\limits_{r>0}\frac{1}{\varphi(r)}\int_{\mathbb{R}^n}\frac{\Phi(|y|)}{|y|^{n}}\Big(\frac{1}{\nu_nr^{n}}\int_{0}^{r}\Big(\int_{{S}^{n-1}}|f(\frac{\rho\theta}{|y|})|^{\tilde{p}_1}d\sigma(\theta)\Big)^{p_1/\tilde{p}_1}\rho^{n-1}d\rho\Big)^{1/p_1}\\
&\quad\times\Big(\frac{1}{\nu_nr^{n}}\int_{0}^{r}\Big(\int_{{S}^{n-1}}|b_{|y|^{-1}B}-b(\frac{\rho\theta}{|y|})|^{\tilde{p}_2}d\sigma(\theta)\Big)^{p_2/\tilde{p}_2}\rho^{n-1}d\rho\Big)^{1/p_2}dy\\
&\leq\sup\limits_{r>0}\frac{1}{\varphi(r)}\int_{\mathbb{R}^n}\frac{\Phi(|y|)}{|y|^{n}}\varphi(r/|y|)dy\|f\|_{{\mathcal{\varphi}L^{p_1}_{\rm rad}L^{\tilde{p}_1}_{\rm ang}(\mathbb R^n)}}\|b\|_{{CMOL^{p_2}_{\rm rad}L^{\tilde{p}_2}_{\rm ang}(\mathbb R^n)}}.
\end{align*}
Combing the above estimates $I,II_1,II_2$ and $III$
leads to our desired conclusion and then completes the proof of Theorem 5.2.
\end{proof}

\noindent\textbf{Theorem 5.3} Let $1\leq p, \tilde{p}<\infty$, if $\Phi$ is a radial function and
$$\mathcal{K}_1=\omega_{n-1}\sup\limits_{r>0}\frac{1}{\varphi(r)}\int_{0}^{\infty}\frac{|\Phi( t)|}{t}\varphi(r/t)dt<\infty.$$
Then
$$\|H_{\Phi}f\|_{{\mathcal{\varphi}L^{p}_{\rm rad}L^{\tilde{p}}_{\rm ang}(\mathbb R^n)}}\leq\mathcal{K}_1\|f\|_{{\mathcal{\varphi}L^{p}_{\rm rad}L^{\tilde{p}}_{\rm ang}(\mathbb R^n)}}.$$
\begin{proof} The Minkowski inequality gives rise to
\begin{align*}
&\|H_{\Phi}f\|_{{\mathcal{\varphi}L^{p}_{\rm rad}L^{\tilde{p}}_{\rm ang}(\mathbb R^n)}}\\
&=\sup\limits_{r>0}\frac{1}{\varphi(r)}\Big(\frac{1}{\nu_nr^{n}}\int_{0}^{r}\Big(\int_{{S}^{n-1}}|\int_{\mathbb{R}^n}\frac{\Phi(\rho\theta|y|^{-1})}{|y|^{n}}f(y)dy|^{\tilde{p}}d\sigma(\theta)\Big)^{p/\tilde{p}}\rho^{n-1}d\rho\Big)^{1/p}\\
&=\omega_{n-1}^{\frac{1}{\tilde{p}}}\sup\limits_{r>0}\frac{1}{\varphi(r)}\Big(\frac{1}{\nu_nr^{n}}\int_{0}^{r}\Big|\int_{0}^{\infty}\int_{S^{n-1}}\frac{\Phi(t)}{t}f(\rho t^{-1}y')d\sigma(y')dt\Big|^{p}\rho^{n-1}d\rho\Big)^{1/p}\\
&\leq\omega_{n-1}^{\frac{1}{\tilde{p}}}\sup\limits_{r>0}\frac{1}{\varphi(r)}\int_{0}^{\infty}\Big(\int_{0}^{r}\frac{1}{\nu_nr^{n}}\Big|\int_{{S}^{n-1}}\frac{\Phi(t)}{t}f(\rho t^{-1}y')d\sigma(y')\Big|^{p}\rho^{n-1}d\rho\Big)^{1/p}dt\\
&\leq\omega_{n-1}^{\frac{1}{\tilde{p}}}\sup\limits_{r>0}\frac{1}{\varphi(r)}\int_{0}^{\infty}\frac{|\Phi( t)|}{t}\Big(\frac{1}{\nu_nr^{n}}\int_{0}^{r/t}\Big|\int_{{S}^{n-1}}f(sy')d\sigma(y')\Big|^{p}s^{n-1}ds\Big)^{1/p}t^{n/p}dt\\
&\leq\omega_{n-1}\sup\limits_{r>0}\frac{1}{\varphi(r)}\int_{0}^{\infty}\frac{|\Phi( t)|}{t}\varphi(r/t)dt\|f\|_{{\mathcal{\varphi}L^{p}_{\rm rad}L^{\tilde{p}}_{\rm ang}(\mathbb R^n)}}.
\end{align*}
This implies Theorem 5.3.
\end{proof}
\noindent\textbf{Theorem 5.4}  Let $1\leq p, p_i, \tilde{p}<\infty$,  $\frac{1}{p}=\frac{1}{p_1}+\frac{1}{p_2}$, $\frac{1}{\tilde{p}}=\frac{1}{\tilde{p}_1}+\frac{1}{\tilde{p}_2}$, if  $\Phi$ is a radial function and
\begin{align*}
\mathcal{K}_1=\omega_{n-1}\sup\limits_{r>0}\frac{1}{\varphi(r)}\int_{0}^{\infty}\frac{|\Phi( t)|}{t}\varphi(r/t)dt<\infty.
\end{align*}
Then
$$
\|H^b_{\Phi}f\|_{{\mathcal{\varphi}L^{p}_{\rm rad}L^{\tilde{p}}_{\rm ang}(\mathbb R^n)}}\lesssim \mathcal{K}_1\|b\|_{{CMOL^{p_2}_{\rm rad}L^{\tilde{p}_2}_{\rm ang}(\mathbb R^n)}}\|f\|_{{\mathcal{\varphi}L^{p_1}_{\rm rad}L^{\tilde{p}_1}_{\rm ang}(\mathbb R^n)}}.
$$
\begin{proof} Note that
\begin{align*}
&\|H^b_{\Phi}f\|_{{\mathcal{\varphi}L^{p}_{\rm rad}L^{\tilde{p}}_{\rm ang}(\mathbb R^n)}}\\
&=\sup\limits_{r>0}\frac{1}{\varphi(r)}\Big(\frac{1}{\nu_nr^{n}}\int_{0}^{r}\Big(\int_{{S}^{n-1}}|\int_{\mathbb{R}^n}\frac{\Phi(\rho\theta|y|^{-1})}{|y|^{n}}f(y)(b(\rho\theta)-b(y))dy|^{\tilde{p}}d\sigma(\theta)\Big)^{p/\tilde{p}}\rho^{n-1}d\rho\Big)^{1/p}\\
&\leq\sup\limits_{r>0}\frac{1}{\varphi(r)}\Big(\frac{1}{\nu_nr^{n}}\int_{0}^{r}\Big(\int_{{S}^{n-1}}|\int_{\mathbb{R}^n}\frac{\Phi(\rho\theta|y|^{-1})}{|y|^{n}}f(\frac{\rho\theta}{|y|})(b(\rho\theta)-b_B)dy|^{\tilde{p}}d\sigma(\theta)\Big)^{p/\tilde{p}}\rho^{n-1}d\rho\Big)^{1/p}\\
&\quad+\sup\limits_{r>0}\frac{1}{\varphi(r)}\Big(\frac{1}{\nu_nr^{n}}\int_{0}^{r}\Big(\int_{{S}^{n-1}}|\int_{\mathbb{R}^n}\frac{\Phi(\rho\theta|y|^{-1})}{|y|^{n}}f(\frac{\rho\theta}{|y|})(b_B-b(y))dy|^{\tilde{p}}d\sigma(\theta)\Big)^{p/\tilde{p}}\rho^{n-1}d\rho\Big)^{1/p}\\
&:=J_1+J_2.
\end{align*}
For $J_1$, the  H\"{o}lder's inequality deduces that
\begin{align*}
&J_1\leq\sup\limits_{r>0}\frac{1}{\varphi(r)}\Big(\frac{1}{\nu_nr^{n}}\int_{0}^{r}\Big(\int_{{S}^{n-1}}|(b(\rho\theta)-b_B)\int_{\mathbb{R}^n}\frac{\Phi(\rho\theta|y|^{-1})}{|y|^{n}}f(\frac{\rho\theta}{|y|})dy|^{\tilde{p}}d\sigma(\theta)\Big)^{p/\tilde{p}}\rho^{n-1}d\rho\Big)^{1/p}\\\
&\leq\sup\limits_{r>0}\frac{1}{\varphi(r)}\Big(\frac{1}{\nu_nr^{n}}\int_{0}^{r}\Big(\int_{{S}^{n-1}}|\int_{\mathbb{R}^n}\frac{\Phi(\rho\theta|y|^{-1})}{|y|^{n}}f(\frac{\rho\theta}{|y|})dy|^{\tilde{p}_1}d\sigma(\theta)\Big)^{1/\tilde{p}_1}\\
&\quad\times\Big(\int_{{S}^{n-1}}|b(\rho\theta)-b_B|^{\tilde{p}_2}d\sigma(\theta)\Big)^{1/\tilde{p}_2}\Big)^p\rho^{n-1}d\rho\Big)^{1/p}\\
&\leq\sup\limits_{r>0}\frac{1}{\varphi(r)}\Big(\frac{1}{\nu_nr^{n}}\int_{0}^{r}\Big(\int_{{S}^{n-1}}|\int_{\mathbb{R}^n}\frac{\Phi(\rho\theta|y|^{-1})}{|y|^{n}}f(\frac{\rho\theta}{|y|})dy|^{\tilde{p}_1}d\sigma(\theta)\Big)^{p_1/\tilde{p}_1}\rho^{n-1}d\rho\Big)^{1/p_1}\\
&\quad\times\Big(\frac{1}{\nu_nr^{n}}\int_{0}^{r}\Big(\int_{{S}^{n-1}}|b(\rho\theta)-b_B|^{\tilde{p}_2}d\sigma(\theta)\Big)^{p_2/\tilde{p}_2}\rho^{n-1}d\rho\Big)^{1/p_2}dy\\
&\leq\|H_{\Phi}f\|_{{\mathcal{\varphi}L^{p_1}_{\rm rad}L^{\tilde{p}_1}_{\rm ang}(\mathbb R^n)}}\|b\|_{{CMOL^{p_2}_{\rm rad}L^{\tilde{p}_2}_{\rm ang}(\mathbb R^n)}}\\
&\leq\mathcal{K}_1\|b\|_{{CMOL^{p_2}_{\rm rad}L^{\tilde{p}_2}_{\rm ang}(\mathbb R^n)}}\|f\|_{{\mathcal{\varphi}L^{p_1}_{\rm rad}L^{\tilde{p}_1}_{\rm ang}(\mathbb R^n)}}.
\end{align*}
For $J_2$, we also have
\begin{align*}
J_2&\leq\|H_{\Phi}(b_B-f)\|_{{\mathcal{\varphi}L^{p}_{\rm rad}L^{\tilde{p}}_{\rm ang}(\mathbb R^n)}}\\
&\leq\mathcal{K}_1\|b_B-f\|_{{\mathcal{\varphi}L^{p}_{\rm rad}L^{\tilde{p}}_{\rm ang}(\mathbb R^n)}}\\
&\leq\mathcal{K}_1\|b\|_{{CMOL^{p_2}_{\rm rad}L^{\tilde{p}_2}_{\rm ang}(\mathbb R^n)}}\|f\|_{{\mathcal{\varphi}L^{p_1}_{\rm rad}L^{\tilde{p}_1}_{\rm ang}(\mathbb R^n)}}.
\end{align*}
Therefore, Theorem 5.3 is finished.
\end{proof}
\section{Boundedness of commutators of  Hausdorff operators on mixed radial-angular Herz type  spaces}

In this section, we will establish the boundedness of commutators of  Hausdorff operators
\begin{align*}
H_{\Phi}^bf(x)&=\int_{\mathbb{R}^n}\frac{\Phi(x|y|^{-1})}{|y|^{n}}f(y)(b(x)-b(y))dy\\
&=\int_{|y|<|x|}\frac{\Phi(x|y|^{-1})}{|y|^{n}}f(y)(b(x)-b(y))dy+\int_{|y|\geq|x|}\frac{\Phi(x|y|^{-1})}{|y|^{n}}f(y)(b(x)-b(y))dy\\
&:=H_{\Phi,1}^bf(x)+H_{\Phi,2}^bf(x).
\end{align*}
on mixed radial-angular Herz type  spaces.

\noindent\textbf{Theorem 6.1} Let $0<q<\infty$, $1<p,\tilde{p}, \tilde{p}_1, \tilde{p}_2<\infty$, $1<r<p$, $\frac{1}{\tilde{p}}=\frac{1}{\tilde{p}_1}+\frac{1}{\tilde{p}_2}$ and $b\in{CMOL^{\max(p,\frac{pr}{p-r})}_{\rm rad}L^{\tilde{p}_1}_{\rm ang}(\mathbb R^n)}$. Then

(a) if $\alpha<n(1/r-1/p)$, $\Phi$ is a radial function and $\Phi$ satisfies
$$\mathcal{A}=\Big(\int_{0}^{\infty}|\Phi(t)|^{r'}\frac{dt}{t^{1-n(r'-1)}}\Big)^{1/r'}<\infty,$$
then $H_{\Phi,1}^b$ is bounded  ${{\dot{K}}^{\alpha,q}_{L^{p}_{\rm rad}L^{\tilde{p}_2}_{\rm ang}}\left(\mathbb R^n\right)}$ to ${{\dot{K}}^{\alpha,q}_{L^{p}_{\rm rad}L^{\tilde{p}}_{\rm ang}}\left(\mathbb R^n\right)}$ .

(b) if $\alpha>-n(1/p-1/r')$, $\Phi$ is a radial function and $\Phi$ satisfies
$$\mathcal{B}=\Big(\int_{0}^{\infty}|\Phi(t)|^{r'}\frac{dt}{t^{1-n}}\Big)^{1/r'}<\infty,$$
then $H_{\Phi,2}^b$ is bounded from  ${{\dot{K}}^{\alpha,q}_{L^{p}_{\rm rad}L^{\tilde{p}_2}_{\rm ang}}\left(\mathbb R^n\right)}$ to ${{\dot{K}}^{\alpha,q}_{L^{p}_{\rm rad}L^{\tilde{p}}_{\rm ang}}\left(\mathbb R^n\right)}$ .

(c) if $-n(1/p-1/r')<\alpha<n(1/r-1/p)$, $\Phi$ is a radial function and $\Phi$ satisfies $\mathcal{A}<\infty$ and $\mathcal{B}<\infty$,
then $H_{\Phi}^b$ is bounded  ${{\dot{K}}^{\alpha,q}_{L^{p}_{\rm rad}L^{\tilde{p}_2}_{\rm ang}}\left(\mathbb R^n\right)}$ to ${{\dot{K}}^{\alpha,q}_{L^{p}_{\rm rad}L^{\tilde{p}}_{\rm ang}}\left(\mathbb R^n\right)}$.

We also generalize the above results to the following.

\noindent\textbf{Theorem 6.2} Let $\lambda\geq0$, $0<q<\infty$, $1<p,\tilde{p}, \tilde{p}_1, \tilde{p}_2<\infty$, $1<r<p$, $\frac{1}{\tilde{p}}=\frac{1}{\tilde{p}_1}+\frac{1}{\tilde{p}_2}$ and $b\in{CMOL^{\max(p,\frac{pr}{p-r})}_{\rm rad}L^{\tilde{p}_1}_{\rm ang}(\mathbb R^n)}$. Then

(a) if $\alpha<\lambda+n(1/r-1/p)$, $\Phi$ is a radial function and $\Phi$ satisfies
$$\mathcal{A}=\Big(\int_{0}^{\infty}|\Phi(t)|^{r'}\frac{dt}{t^{1-n(r'-1)}}\Big)^{1/r'}<\infty,$$
then $H_{\Phi,1}^b$ is bounded from ${M\dot{K}}^{\alpha,q,\lambda}_{L^{p}_{\rm rad}L^{\tilde{p}_2}_{\rm ang}}(\mathbb R^{n})$ to ${M\dot{K}}^{\alpha,q,\lambda}_{L^{p}_{\rm rad}L^{\tilde{p}}_{\rm ang}}(\mathbb R^{n})$.

(b) if $\alpha>\lambda-n(1/p-1/r')$, $\Phi$ is a radial function and $\Phi$ satisfies
$$\mathcal{B}=\Big(\int_{0}^{\infty}|\Phi(t)|^{r'}\frac{dt}{t^{1-n}}\Big)^{1/r'}<\infty,$$
then $H_{\Phi,2}^b$ is bounded from ${M\dot{K}}^{\alpha,q,\lambda}_{L^{p}_{\rm rad}L^{\tilde{p}_2}_{\rm ang}}(\mathbb R^{n})$ to ${M\dot{K}}^{\alpha,q,\lambda}_{L^{p}_{\rm rad}L^{\tilde{p}}_{\rm ang}}(\mathbb R^{n})$.

(c) if $\lambda-n(1/p-1/r')<\alpha<\lambda+n(1/r-1/p)$, $\Phi$ is a radial function and $\Phi$ satisfies $\mathcal{A}<\infty$ and $\mathcal{B}<\infty$,
then $H_{\Phi}^b$ is bounded from ${M\dot{K}}^{\alpha,q,\lambda}_{L^{p}_{\rm rad}L^{\tilde{p}_2}_{\rm ang}}(\mathbb R^{n})$ to ${M\dot{K}}^{\alpha,q,\lambda}_{L^{p}_{\rm rad}L^{\tilde{p}}_{\rm ang}}(\mathbb R^{n})$.

In fact, the proofs of Theorem 6.1 and Theorem 6.2 are similar, so we only need provide the proof of  Theorem 6.1.
\begin{proof} Note that
\begin{align*}
&\|H_{\Phi,1}^bf\chi_k\|^p_{L^{p}_{\rm rad}L^{\tilde{p}}_{\rm ang}(\mathbb R^{n})}\\
&=\int_{2^{k-1}}^{2^k}\Big(\int_{{S}^{n-1}}\Big|\int_{|y|<\rho}\frac{\Phi(\rho|y|^{-1})}{|y|^{n}}f(y)(b(\rho\theta)-b(y))dy\Big|^{\tilde{p}}d\sigma(\theta)\Big)^{p/\tilde{p}}\rho^{n-1}d\rho\\
&\lesssim\int_{2^{k-1}}^{2^k}\Big(\int_{{S}^{n-1}}\Big|\int_{|y|<\rho}\frac{\Phi(\rho|y|^{-1})}{|y|^{n}}f(y)(b(\rho\theta)-b_{B_k})dy\Big|^{\tilde{p}}d\sigma(\theta)\Big)^{p/\tilde{p}}\rho^{n-1}d\rho\\
&\quad+\int_{2^{k-1}}^{2^k}\Big|\int_{|y|<\rho}\frac{\Phi(\rho|y|^{-1})}{|y|^{n}}f(y)(b(y)-b_{B_k})dy\Big|^{p}\rho^{n-1}d\rho\\
&:=A_1+A_2.
\end{align*}
For $A_1$, the Minkowski inequality and  H\"{o}lder's inequality deduce that
\begin{align*}
&\int_{2^{k-1}}^{2^k}\Big|\int_{|y|<\rho}\frac{\Phi(\rho|y|^{-1})}{|y|^{n}}f(y)dy\Big|^p\Big(\int_{{S}^{n-1}}|b(\rho\theta)-b_{B_k}|^{\tilde{p}}d\sigma(\theta)\Big)^{p/\tilde{p}}\rho^{n-1}d\rho\\
&=\int_{2^{k-1}}^{2^k}\Big|\sum\limits_{i=-\infty}^{k}\int_{2^{i-1}}^{2^i}\frac{\Phi(\rho t^{-1})}{t}\int_{{S}^{n-1}}f(ty')d\sigma(y')dt\Big|^p\Big(\int_{{S}^{n-1}}|b(\rho\theta)-b_{B_k}|^{\tilde{p}}d\sigma(\theta)\Big)^{p/\tilde{p}}\rho^{n-1}d\rho\\
&\lesssim\int_{2^{k-1}}^{2^k}\Big(\sum\limits_{i=-\infty}^{k}\int_{2^{i-1}}^{2^i}\frac{|\Phi(\rho t^{-1})|}{t}\Big(\int_{{S}^{n-1}}|f(ty')|^{\tilde{p}}d\sigma(y')\Big)^{1/\tilde{p}}dt\Big)^p\\
&\quad\times\Big(\int_{{S}^{n-1}}|b(\rho\theta)-b_{B_k}|^{\tilde{p}}d\sigma(\theta)\Big)^{p/\tilde{p}}\rho^{n-1}d\rho\\
&\lesssim\int_{2^{k-1}}^{2^k}\sum\limits_{i=-\infty}^{k}\Big(\int_{2^{i-1}}^{2^i}\frac{|\Phi(\rho t^{-1})|^{r'}}{t^{nr'}}t^{n-1}dt\Big)^{p/r'}\Big(\int_{2^{i-1}}^{2^i}\Big(\int_{{S}^{n-1}}|f(ty')|^{\tilde{p}}d\sigma(y')\Big)^{r/\tilde{p}}t^{n-1}dt\Big)^{p/r}\\
&\quad\times\Big(\int_{{S}^{n-1}}|b(\rho\theta)-b_{B_k}|^{\tilde{p}}d\sigma(\theta)\Big)^{p/\tilde{p}}\rho^{n-1}d\rho\\
&\lesssim\int_{2^{k-1}}^{2^k}\Big(\sum\limits_{i=-\infty}^{k}\Big(\int_{\rho/2^{i}}^{\rho/2^{i-1}}|\Phi(t)|^{r'}t^{-1+n(r'-1)}\rho^{n(1-r')}dt\Big)^{p/r'}\|f\chi_i\|^p_{{L^{p}_{\rm rad}L^{\tilde{p}}_{\rm ang}(\mathbb R^n)}}\Big)\\
&\quad\times 2^{inp(1/r-1/p)}\Big(\int_{{S}^{n-1}}|b(\rho\theta)-b_{B_k}|^{\tilde{p}}d\sigma(\theta)\Big)^{p/\tilde{p}}\rho^{n-1}d\rho\\
&\lesssim\mathcal{A}^p\Big(\sum\limits_{i=-\infty}^{k}2^{in(1/r-1/p)}2^{-kn/r}\|f\chi_i\|_{{L^{p}_{\rm rad}L^{\tilde{p}}_{\rm ang}(\mathbb R^n)}}\Big)^p 2^{kn}\|b\|^p_{{CMOL^{p}_{\rm rad}L^{\tilde{p}}_{\rm ang}(\mathbb R^n)}}\\
&\lesssim\mathcal{A}^p\Big(\sum\limits_{i=-\infty}^{k}2^{(i-k)n(1/r-1/p)}\|f\chi_i\|_{{L^{p}_{\rm rad}L^{\tilde{p}}_{\rm ang}(\mathbb R^n)}}\Big)^p \|b\|^p_{{CMOL^{p}_{\rm rad}L^{\tilde{p}}_{\rm ang}(\mathbb R^n)}}.
\end{align*}
For $A_2$, we divide it into two parts as follows:
\begin{align*}
A_2&\lesssim\int_{2^{k-1}}^{2^k}\Big|\sum\limits_{i=-\infty}^{k}\int_{2^{i-1}}^{2^i}\frac{\Phi(\rho|y|^{-1})}{|y|^{n}}f(y)(b(y)-b_{B_i})dy\Big|^{p}\rho^{n-1}d\rho\\
&\quad+\int_{2^{k-1}}^{2^k}\Big|\sum\limits_{i=-\infty}^{k}\int_{2^{i-1}}^{2^i}\frac{\Phi(\rho|y|^{-1})}{|y|^{n}}f(y)(k-i)\|b\|_{{CMOL^{1}_{\rm rad}L^{\tilde{p}}_{\rm ang}(\mathbb R^n)}}dy\Big|^{p}\rho^{n-1}d\rho\\
&:=A_{21}+A_{22}.
\end{align*}
For $A_{21}$, the Minkowski inequality and  H\"{o}lder's inequality give rise to
\begin{align*}
&A_{21}\leq\int_{2^{k-1}}^{2^k}\Big(\sum\limits_{i=-\infty}^{k}\Big(\int_{2^{i-1}}^{2^i}\frac{|\Phi(\rho t^{-1})|^{r'}}{t^{nr'}}t^{n-1}dt\Big)^{1/r'}\\
&\times\Big(\int_{2^{i-1}}^{2^i}\Big(\int_{{S}^{n-1}}|f(ty')(b(ty')-b_{B_i})|^{\tilde{p}}d\sigma(y')\Big)^{r/\tilde{p}}t^{n-1}dt\Big)^{1/r}\Big)^{p}\rho^{n-1}d\rho\\
&\lesssim\mathcal{A}^p\int_{2^{k-1}}^{2^k}\Big[\sum\limits_{i=-\infty}^{k}\rho^{-n/r}\Big(\int_{2^{i-1}}^{2^i}\int_{{S}^{n-1}}|f(ty')(b(ty')-b_{B_i})|^{\tilde{p}}d\sigma(y')\Big)^{r/\tilde{p}}t^{n-1}dt\Big)^{1/r}\Big]^p\rho^{n-1}d\rho\\
&\lesssim\mathcal{A}^p\int_{2^{k-1}}^{2^k}\Big[\sum\limits_{i=-\infty}^{k}\rho^{-n/r}\Big(\int_{2^{i-1}}^{2^i}\Big(\int_{{S}^{n-1}}|b(ty')-b_{B_i}|^{\tilde{p}_1}d\sigma(y')\Big)^{rp/(p-r)\tilde{p}_1}t^{n-1}dt\Big)^{p-r/pr}\\
&\quad\times\Big(\int_{2^{i-1}}^{2^i}\Big(\int_{{S}^{n-1}}|f(ty')|^{\tilde{p}_2}d\sigma(y')\Big)^{p/\tilde{p}_2}t^{n-1}dt\Big)^{1/p}\Big]^p\rho^{n-1}d\rho\\
&\lesssim\mathcal{A}^p\Big(\sum\limits_{i=-\infty}^{k}2^{kn(1/p-1/r)}\|b\|_{{CMOL^{pr/p-r}_{\rm rad}L^{\tilde{p}_1}_{\rm ang}(\mathbb R^n)}}2^{in(p-r/pr)}\|f\chi_i\|_{{L^{p}_{\rm rad}L^{\tilde{p}_2}_{\rm ang}(\mathbb R^n)}}\Big)^p\\
&\lesssim\mathcal{A}^p\|b\|^p_{{CMOL^{pr/p-r}_{\rm rad}L^{\tilde{p}_1}_{\rm ang}(\mathbb R^n)}}\Big(\sum\limits_{i=-\infty}^{k}2^{(i-k)n(1/p-1/r)}\|f\chi_i\|_{{L^{p}_{\rm rad}L^{\tilde{p}_2}_{\rm ang}(\mathbb R^n)}}\Big)^p.
\end{align*}
For $A_{22}$, we also have
\begin{align*}
&\int_{2^{k-1}}^{2^k}\Big|\sum\limits_{i=-\infty}^{k}\int_{2^{i-1}}^{2^i}\frac{\Phi(\rho|y|^{-1})}{|y|^{n}}f(y)(k-i)\|b\|_{{CMOL^{1}_{\rm rad}L^{\tilde{p}}_{\rm ang}(\mathbb R^n)}}dy\Big|^{p}\rho^{n-1}d\rho\\
&\lesssim\|b\|^p_{{CMOL^{1}_{\rm rad}L^{\tilde{p}}_{\rm ang}(\mathbb R^n)}}\Big(\sum\limits_{i=-\infty}^{k}(k-i)2^{(i-k)n(1/r-1/p)}\|f\chi_i\|_{{L^{p}_{\rm rad}L^{\tilde{p}}_{\rm ang}(\mathbb R^n)}}\Big)^p.
\end{align*}
Therefore, we have
\begin{align*}
&\|H_{\Phi,1}^bf\chi_k\|_{{\dot{K}}^{\alpha,q}_{L^{p}_{\rm rad}L^{\tilde{p}}_{\rm ang}}\left(\mathbb R^n\right)}\\
&\lesssim\mathcal{A}\|b\|_{{CMOL^{p}_{\rm rad}L^{\tilde{p}}_{\rm ang}(\mathbb R^n)}}\Big\{\sum\limits_{k=-\infty}^{\infty}2^{k\alpha q} \Big(\sum\limits_{i=-\infty}^{k}2^{(i-k)n(1/r-1/p)}\|f\chi_i\|_{{L^{p}_{\rm rad}L^{\tilde{p}}_{\rm ang}(\mathbb R^n)}}\Big)^q\Big\}^{1/q}\\
&\quad+\mathcal{A}\|b\|_{{CMOL^{pr/p-r}_{\rm rad}L^{\tilde{p}_1}_{\rm ang}(\mathbb R^n)}}\Big\{\sum\limits_{k=-\infty}^{\infty}2^{k\alpha q} \Big(\sum\limits_{i=-\infty}^{k}2^{(i-k)n(1/r-1/p)}\|f\chi_i\|_{{L^{p}_{\rm rad}L^{\tilde{p}_2}_{\rm ang}(\mathbb R^n)}}\Big)^q\Big\}^{1/q}\\
&\quad+\mathcal{A}\|b\|_{{CMOL^{1}_{\rm rad}L^{\tilde{p}}_{\rm ang}(\mathbb R^n)}}\Big\{\sum\limits_{k=-\infty}^{\infty}2^{k\alpha q} \Big(\sum\limits_{i=-\infty}^{k}(k-i)2^{(i-k)n(1/r-1/p)}\|f\chi_i\|_{{L^{p}_{\rm rad}L^{\tilde{p}}_{\rm ang}(\mathbb R^n)}}\Big)^q\Big\}^{1/q}\\
&\lesssim\mathcal{A}\|b\|_{CMOL^{\max(p,\frac{pr}{p-r})}_{\rm rad}L^{\tilde{p}_1}_{\rm ang}(\mathbb R^n)}\|f\|_{{\dot{K}}^{\alpha,q}_{L^{p}_{\rm rad}L^{\tilde{p}_2}_{\rm ang}}\left(\mathbb R^n\right)}.
\end{align*}
Theorem 5.3 (a) is finished.

The proof of (b) is similar to the proof of (a), and combining (a) with (b) leads to (c). This finishes the proof of Theorem 5.3.
\end{proof}
\hspace*{\fill} \\
\noindent
\textbf{Acknowledgments}

The first author was supported by the Doctoral Scientific Research Foundation of Northwest Normal University (No. 202203101202), the Young Teachers Scientific Research Ability Promotion Project of Northwest Normal University (No. NWNU-LKQN2023-15)
and the Key Laboratory of Computing and Stochastic Mathematics (Ministry of Education), School of Mathematics and Statistics, Hunan Normal University (No. 202402).








\end{document}